\pdfoutput=1
\documentclass[reqno]{amsart}

\usepackage{graphicx}
\usepackage{amssymb}
\usepackage{amsmath}
\usepackage{hyperref}
\newcommand{\abs}[1]{\left|#1\right|}

\newcommand{\vnormf}[1]{||#1||}

\newcommand{\sign}[1]{\mbox{sign}#1}
\newcommand{\Z}{\mathbb{Z}}

\newcommand{\R}{\mathbb{R}}

\newcommand{\onewidth}{.5}
\newcommand{\onewidthh}{.6}

\newcommand{\clpg}{\clearpage}
\newcommand{\lbreak}{\\}
\newcommand{\fourteenwidth}{.7}
\newcommand{\finonewidth}{1}
\newcommand{\fintwowidth}{1}
\newcommand{\finthreewidth}{1}
\newcommand{\finfourwidth}{1}
\long\def\symbolfootnote[#1]#2{\begingroup%
\def\thefootnote{\fnsymbol{footnote}}\footnote[#1]{#2}\endgroup}

\begin{document}

\title{Orthogonal Polynomials with Respect to Self-Similar Measures}
\author{Steven M. Heilman$^{1}$}
\author{Philip Owrutsky$^{2}$}
\author{Robert S. Strichartz$^{3}$}
\thanks{
Department of Mathematics, Cornell University, Ithaca, NY 14850-4201
\\
Courant Institute of Mathematical Sciences, New York University, New York, NY 10012-1185
\\
\textit{E-mail address}: \texttt{heilman@cims.nyu.edu}
\\
Division of Engineering and Applied Sciences, Harvard University, Cambridge, MA 02138
\\
\textit{E-mail address}: \texttt{owrutsky@fas.harvard.edu}
\\
Department of Mathematics, Cornell University, Ithaca, NY 14850-4201
\\
\textit{E-mail address}: \texttt{str@math.cornell.edu}
\\
$^{1}$: Research supported by the National Science Foundation through the Research Experiences for Undergraduates Program at Cornell.
\\
$^{2}$: Cornell Presidential Research Scholar
\\
$^{3}$: Research supported in part by the National Science Foundation, grant DMS-0652440.
}

\date{\today}

\begin{abstract}
We study experimentally systems of orthogonal polynomials with respect to self-similar measures.  When the support of the measure is a Cantor set, we observe some interesting properties of the polynomials, both on the Cantor set and in the gaps of the Cantor set.  We introduce an effective method to visualize the graph of a function on a Cantor set.  We suggest a new perspective, based on the theory of dynamical systems, for studying families $P_{n}(x)$ of orthogonal functions as functions of $n$ for fixed values of $x$.
\end{abstract}
\maketitle

\section{Introduction}
\label{secintro}

The classical theory of orthogonal polynomials \cite{gautschi04}\cite{szego75} allows you to start the Gram-Schmidt process with virtually any measure.  In the families of polynomials that are usually studied, the measure is either absolutely continuous or discrete, but the general theory allows one to use a singular continuous measure.  In recent years there has been an interest in the case of fractal measures.  See for example \cite{barnsley83}-\cite{barnsley85},\cite{bessis83},\lbreak\cite{mantica96}-\cite{mantica00}. In particular, G. Mantica has developed algorithms to efficiently compute the coefficients of the $3$-term recursion relation, and hence the polynomials, in the case of a self-similar measure \cite{mantica96,mantica00}.  In this work we use these algorithms as a tool to study the polynomials experimentally, looking for interesting patterns and conjectures.  We view this as part of a general program to explore portions of classical analysis related to fractal measures.

Let $\mu$ be a measure on the line, and for simplicity assume that $\mu$ is a probability measure supported on the unit interval.  The related system $\{P_{n}(x)\}$ of orthogonal polynomials is characterized as follows:
\begin{equation}\label{one1}
P_{n}(x)\mbox{ is a polynomial of degree $n$ with a non-negative coefficient of }x^{n},
\end{equation}
\begin{equation}\label{one2}
\int P_{n}(x)P_{m}(x)d\mu(x)=\delta_{n,m}\qquad\mbox{, orthonormality.}
\end{equation}
The general theory implies that there is a $3$-term recursion relation
\begin{equation}\label{one3}
xP_{n}(x)=r_{n+1}P_{n+1}(x)+A_{n}P_{n}(x)+r_{n}P_{n-1},
\end{equation}
where $r_{n},A_{n}>0$, $n\in\Z^{+}$ are determined by the measure $\mu$ and $P_{-1}=0$, $P_{0}=1$.  The coefficients $r_{n},A_{n}$ are also called the entries of the \underline{Jacobi matrix}, an infinite symmetric tri-diagonal matrix $J$ such that $J_{i,i}=A_{i}$ and $J_{i,i+1}=J_{i+1,i}=r_{i+1}$ for $i\in\Z^{\geq0}$.  Note that Eq. (\ref{one3}) allows us to find the polynomials recursively:
\begin{equation}\label{one4}
P_{n+1}(x)=\frac{1}{r_{n+1}}((x-A_{n})P_{n}(x))-\frac{r_{n}}{r_{n+1}}P_{n-1}(x)
\end{equation}
We should point out that while in principle Eq. (\ref{one4}) allows computation of the coefficients of $P_{n}$, this computation could be unstable.  In the case that $A_{n}$ is constant (as we see below when our measure is supported on a Cantor set), it is often easier to compute the coefficients of $P_{n}$ as a polynomials of $(x-A_{n})$.  When the $A_{n}$ vary, we alternatively use Eq. (\ref{one4}) to compute the values $P_{n}(x)$ for specific $x$-values.

A measure $\mu$ on the line is said to be self-similar if there exists an iterated function system (IFS) of contractive similarities $\{F_{i}\}_{i=1}^{N}$ and a set of probability weights $\{p_{i}\}_{i=1}^{N}$ such that
\begin{equation}\label{one5}
\mu(A)=\sum_{i=1}^{N}p_{i}\mu(F_{i}^{-1}A)\mbox{, for any measurable set $A$, or equivalently,}
\end{equation}
\begin{equation}\label{one6}
\int f\,d\mu=\sum_{i=1}^{N}p_{i}\int f\circ F_{i}\,d\mu
\mbox{, for any continuous function }f.
\end{equation}
In this paper we restrict our attention to the family of IFS's with $N=2$ and such that
\begin{equation}\label{one7}
F_{1}(x)=\frac{1}{R}x,\qquad F_{2}(x)=\frac{1}{R}(x-1)+1
\end{equation}
where $R\geq2$ is a parameter.  When $R=2$ and $p_{1}=P_{2}=\frac{1}{2}$ we obtain Lebesgue measure on $[0,1]$, and the corresponding polynomials are essentially the classical Legendre polynomials.  (Actually, the classical Legendre polynomials are orthogonal with respect to Lebesgue measure on $[-1,1]$ and are normalized differently, but the differences just involve rescaling the axes.)  When $R=2$ and $p_{1}\neq p_{2}$, we refer to the corresponding polynomials as \textit{Weighted Legendre Polynomials} (WLP).  The measure $\mu$ is singular but not supported on any proper closed subset of $[0,1]$.  When $R>2$ we will always take $p_{1}=p_{2}=\frac{1}{2}$, and we call the corresponding polynomials \textit{Cantor Legendre Polynomials} (CLP).  The measure $\mu$ is then supported on a Cantor set $C_{R}$ characterized by
\begin{equation}\label{one8}
C_{R}=F_{1}C_{R}\cup F_{2}C_{R}
\end{equation}
The standard Cantor set and Cantor measure correspond to $R=3$.  We refer to the intervals in $[0,1]\setminus C_{R}$ as gaps.  The largest gap is the interval $(\frac{1}{R},1-\frac{1}{R})$, and there are $2^{m}$ gaps of length $\frac{1}{R^{m-1}}(1-\frac{2}{R})$.

The behavior of $P_{n}(x)$ in the CLP case is quite different on the gaps and on the Cantor set $C_{R}$.  In order to visualize the graphs of $P_{n}(x)$ on $C_{R}$ we introduce the distorted Cantor set $\widetilde{C}_{R,\epsilon}$ (also known as the Smith-Volterra-Cantor set or fat Cantor set), obtained by reducing the size of the gaps by a factor of $\epsilon$ (a parameter that we choose).  Note that $\widetilde{C}_{R,\epsilon}$ is still a topological Cantor set, but it has positive Lebesgue measure.  There is an obvious one-to-one correspondence between $C_{R}$ and $\widetilde{C}_{R,\epsilon}$ that identifies regions between corresponding gaps.  We use this identification to graph functions defined on $C_{R}$ against $\widetilde{C}_{R,\epsilon}$.

In section \ref{sectwo} we present data for the entries of the Jacobi matrix.  In the CLP case we note the different behaviors of $r_{n}$ for even and odd $n$, and the occurrence of small values.  In section \ref{secthree} we display graphs of the polynomials.  In the CLP case we show graphs of the restrictions to the Cantor set and to the gaps.

We then discuss various features of the data.  In section \ref{secdir} we discuss the associated Dirichlet kernels.  By using the Christoffel-Darboux formula we are able to relate approximate identity behavior with small values of $r_{n}$.  In section \ref{secapprox} we discuss some approximate equalities relating CLP restricted to the Cantor set.  In particular, $P_{2n}(x)$ is approximately equal to $P_{2n+1}(x)$ on the right half of $C_{R}$, for large $n$.  We also define a ``shuffle'' map that approximately preserves $P_{2n+1}$ when $n$ is a power of $2$.  In section \ref{secclp} we discuss the behavior of $P_{n}(x)$ on the gaps in the CLP case.  On the central gap, for high $n$, $P_{n}(x)$ vaguely approximates either a Gaussian or the derivative of a Gaussian, depending on $n\mod2$.  More precisely, we find that $P_{2n}(x)=c_{2n}e^{-d_{2n}x^{\alpha(x)}}$ where $\alpha(x)=\alpha_{2n}(x)=2$ at $x=1/2$, and $\alpha(x)\approx2$ for $x$ around $1/2$.  On the other gaps, for large enough $n$, the behavior is roughly the same.  In section \ref{secwlp} we discuss the behavior of $P_{n}(x)$ at the points $x=0,1/4,1/2,3/4,1$ in the WLP case, and we contrast the results with the known behavior for Legendre polynomials.  This study leads to the dynamical systems perspective.

Instead of thinking of $P_{n}(x)$ as a function of $x$ for fixed $n$, we look at $P_{n}(x)$ as a function of $n$, for a fixed $x$.  Because of the $3$-term recursion relation, it is more natural to look at vectors $(P_{n}(x),P_{n+1}(x))$ in the plane.  Then there exist $2\times 2$ matrices $M_{n}(x)$ such that
\begin{equation}\label{one9}
\left(\begin{array}{c}P_{n+1}(x)\\P_{n}(x)\end{array}\right)
=M_{n}(x)\left(\begin{array}{c}P_{n}(x)\\P_{n-1}(x)\end{array}\right).
\end{equation}
In fact,
\begin{equation}\label{one10}
M_{n}(x)=
\left(
\begin{array}{cc}
\frac{x-A_{n}}{r_{n+1}} & -\frac{r_{n}}{r_{n+1}}\\
1 & 0
\end{array}
\right).
\end{equation}
So, we are really looking at orbits of a time-dependent linear dynamical system.  We can then define the analogy of the Mandelbrot set as all $x$ for which the orbit is bounded.  It appears that this coincides, more or less, with the support of the measure.  Then for each $x$ in this Mandelbrot set, we can define a Julia set $J(x)$ as the limit set of the orbit.  We display some examples of these Julia sets.  We then extend the investigation of section \ref{secwlp} to generic $x$ for WLP and CLP in section \ref{secdyn}.  We end with a short concluding discussion in Section \ref{secconc}.

This paper should be viewed in the context of a long term effort to understand topics in classical analysis extended to fractal measures.  The following references are just a sampling of this work: \cite{bird03,coletta92,dutkay06,huang01,jorgensen98,kigami01,laba02,lau93,lund98}\lbreak\cite{strichartz90}-\cite{strichartz06}. More data may be found at\lbreak \url{www.math.cornell.edu/~orthopoly}.

\underline{Acknowledgement:} We are grateful to Giorgio Mantica for allowing us to use his codes.

\section{Entries of the Jacobi Matrix}
\label{sectwo}

\begin{figure}[htbp!]
\includegraphics[width=\onewidth\textwidth]{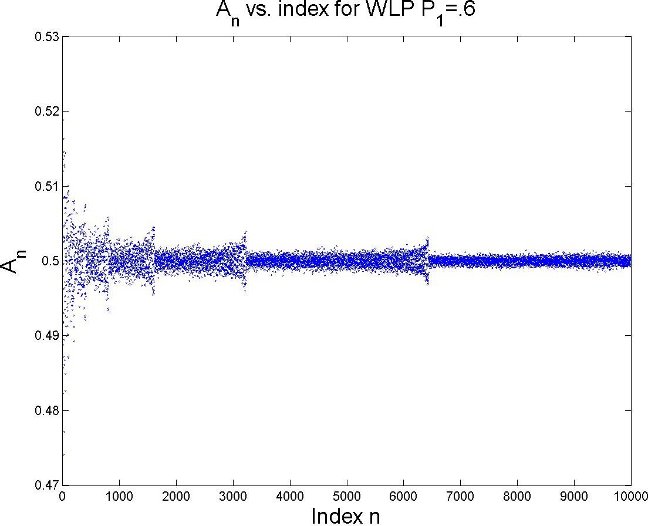}
\caption{Plot of $A_{n}$ vs. index $n$, $1\leq n\leq10000$ for $p_{1}=.6$}
\label{figtwo1}
\end{figure}

\begin{figure}[htbp!]
\includegraphics[width=\onewidth\textwidth]{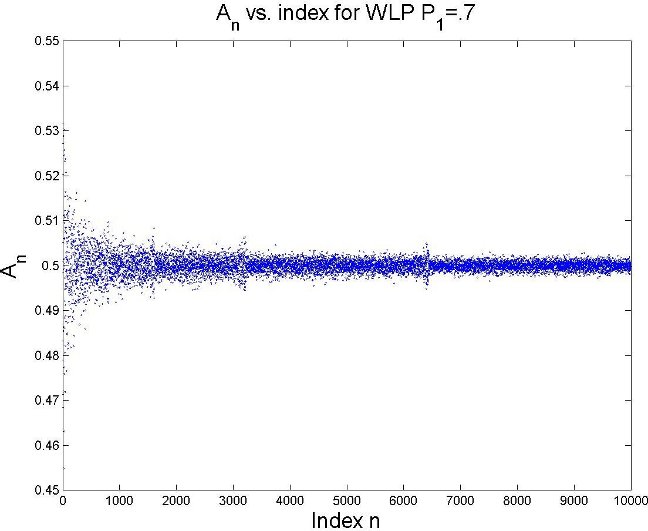}
\caption{Plot of $A_{n}$ vs. index $n$, $1\leq n\leq10000$ for $p_{1}=.7$}
\label{figtwo2}
\end{figure}

\begin{figure}[htbp!]
\includegraphics[width=\onewidth\textwidth]{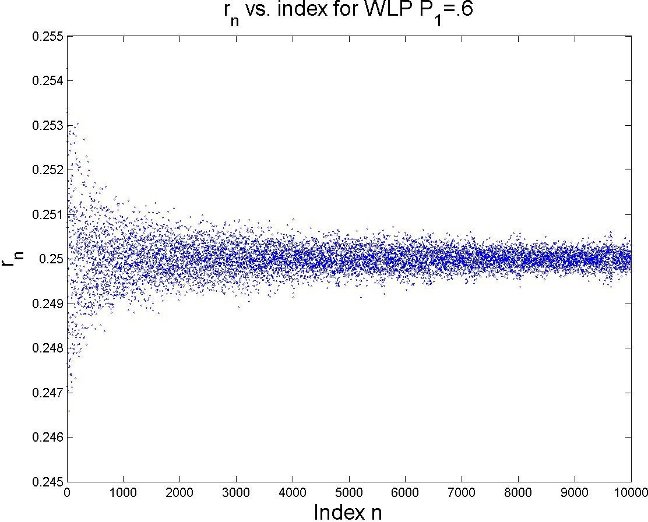}
\caption{Plot of $r_{n}$ vs. index $n$, $1\leq n\leq10000$ for $p_{1}=.6$}
\label{figtwo3}
\end{figure}

\begin{figure}[htbp!]
\includegraphics[width=\onewidth\textwidth]{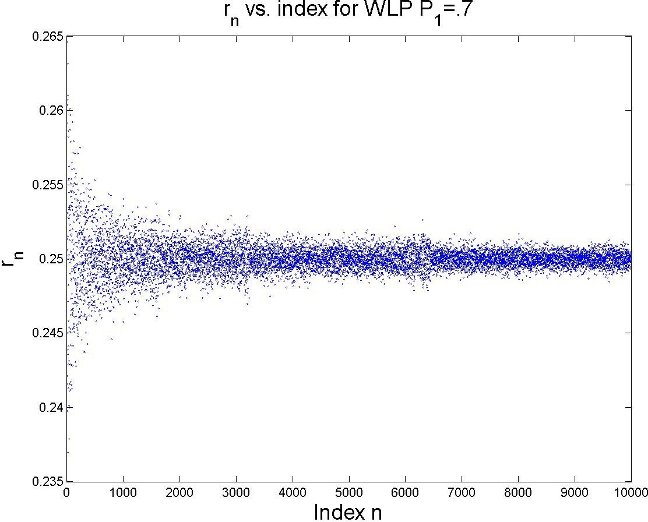}
\caption{Plot of $r_{n}$ vs. index $n$, $1\leq n\leq10000$ for $p_{1}=.7$}
\label{figtwo4}
\end{figure}

The coefficients $r_{n}$ and $A_{n}$ in the $3$-term recursion relation (Eq. (\ref{one3})) determine the polynomials in a rather subtle way.  More work is needed to clarify this relationship.  In this section we report data for our two classes of examples.  In Figs. \ref{figtwo1}-\ref{figtwo4} we graph $A_{n}$ and $r_{n}$ versus $n$, for several choices of $p_{1}$ in the WLP case.  For classical Legendre polynomials we have $A_{n}=1/2$ for all $n$ by symmetry about $x=1/2$, and $\lim_{n\to\infty}r_{n}=1/4$.  The WLP case shows a small but significant difference from this model case.

\begin{figure}[htbp!]
\includegraphics[width=\onewidth\textwidth]{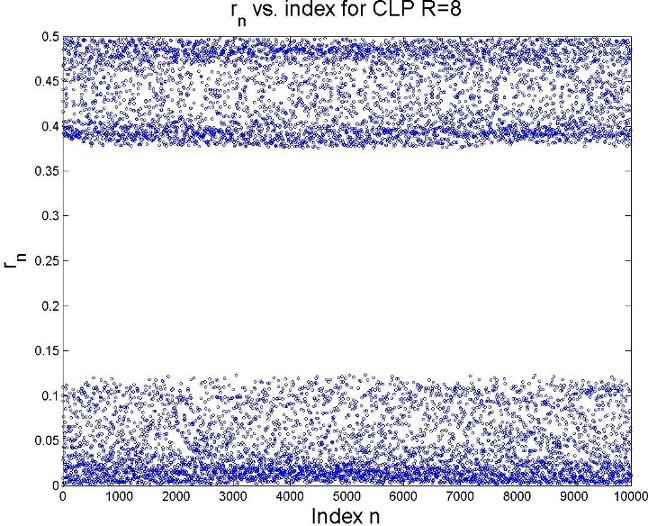}
\caption{Plot of $r_{n}$ vs. index $n$, $1\leq n\leq10000$ for CLP, $R=8$}
\label{figtwo5}
\end{figure}

In the CLP case, all $A_{n}=1/2$ by symmetry about $x=1/2$.  In Fig. \ref{figtwo5} we graph $r_{n}$ versus $n$ for $R=8$.

\underline{Problem 2.1}: What is the nature of the limit set
\begin{equation}\label{two11}
\cap_{N=1}^{\infty}\mbox{cl}\left(\cup_{n=N}^{\infty}\{r_{n}\}\right)?
\end{equation}
Do the measures
\begin{equation}\label{two12}
\frac{1}{m-k}\sum_{n=k+1}^{m}\delta_{r_{n}}
\end{equation}
converge weakly to some fractal measure as $m,k$ tend to infinity in some specific manner?

\begin{figure}[htbp!]
\includegraphics[width=\onewidth\textwidth]{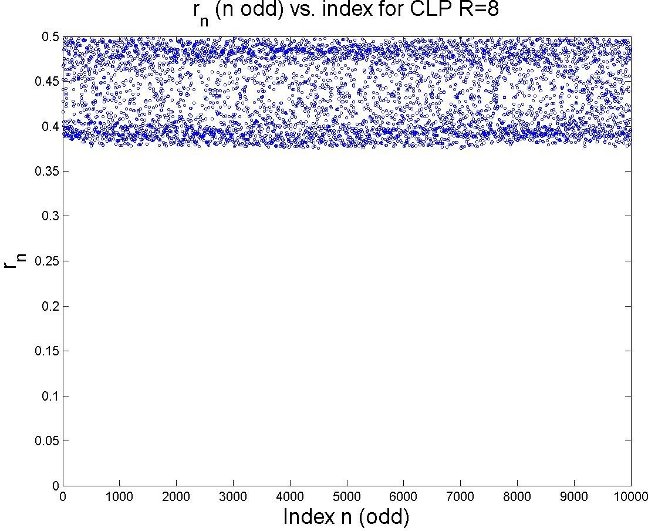}
\caption{Plot of $r_{n}$ vs. index $n$, $1\leq n\leq10000$ for $n$ odd, CLP, $R=8$}
\label{figtwo6}
\end{figure}

\begin{figure}[htbp!]
\includegraphics[width=\onewidth\textwidth]{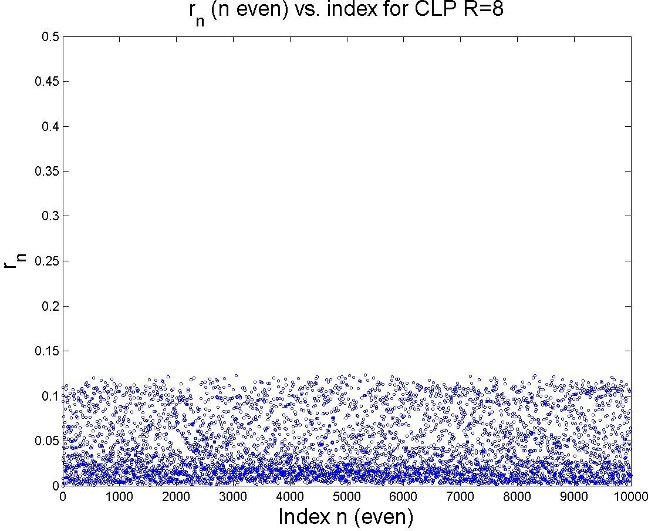}
\caption{Plot of $r_{n}$ vs. index $n$, $1\leq n\leq1000$ for $n$ even, CLP, $R=8$}
\label{figtwo7}
\end{figure}

A striking feature of the data is the different behavior of $r_{n}$ for $n$ even and $n$ odd.  In Figs. \ref{figtwo6} and \ref{figtwo7} we show the same data as Fig. \ref{figtwo5}, separating the even and odd values of $n$.  Another striking observation is that some values of $r_{n}$ for $n$ even are close to zero.

\underline{Problem 2.2}: What is
\begin{equation}\label{two13}
\liminf_{n\to\infty}r_{n}?
\end{equation}
In particular, is the $\liminf$ zero?  What is the sequence of $n$'s along which the $\liminf$ is attained?

As we will see later, having values of $r_{n}$ close to zero has interesting implications.  We could also ask for the $\limsup$, but it is not clear what significance this has.

\underline{Conjecture 2.3}: We always have
\begin{equation}\label{two14}
r_{n}\leq\frac{1}{2}
\end{equation}

\clpg
\section{Graphs of Polynomials}
\label{secthree}

\begin{figure}[htbp!]
\includegraphics[width=\onewidth\textwidth]{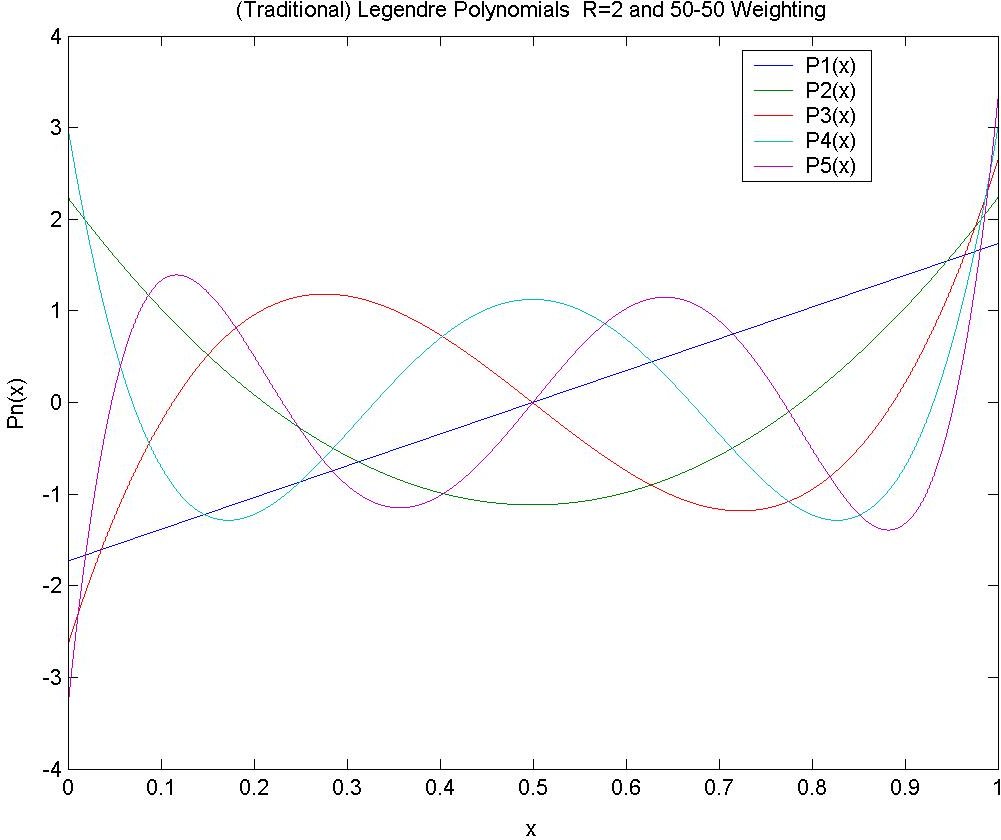}
\caption{Plot of $P_{n}(x)$  on the unit interval with $p_{1}=p_{2}=.5$.  This gives the classical Legendre Polynomials (up to renormalization).}
\label{figthree1}
\end{figure}

\begin{figure}[htbp!]
\includegraphics[width=\onewidth\textwidth]{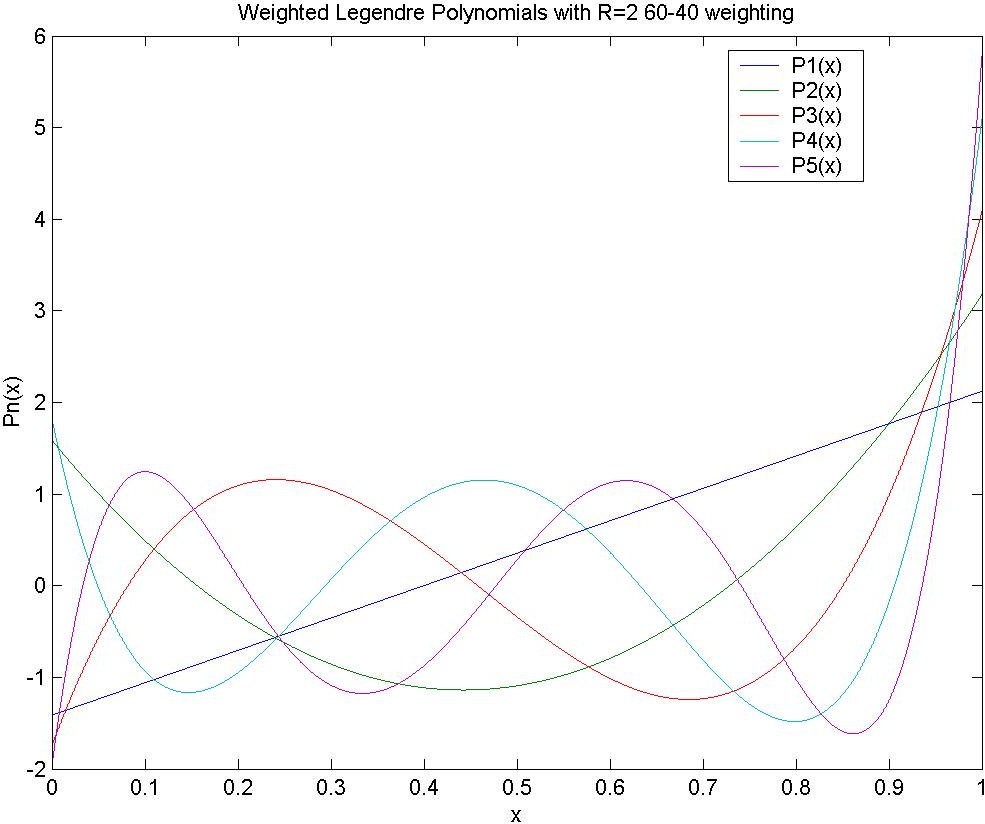}
\caption{Plot of $P_{n}(x)$  on the unit interval with $p_{1}=.6, p_{2}=.4$.  Recall that the $60$ percent weighting is given to the map which contracts towards zero.}
\label{figthree2}
\end{figure}

\begin{figure}[htbp!]
\includegraphics[width=\onewidth\textwidth]{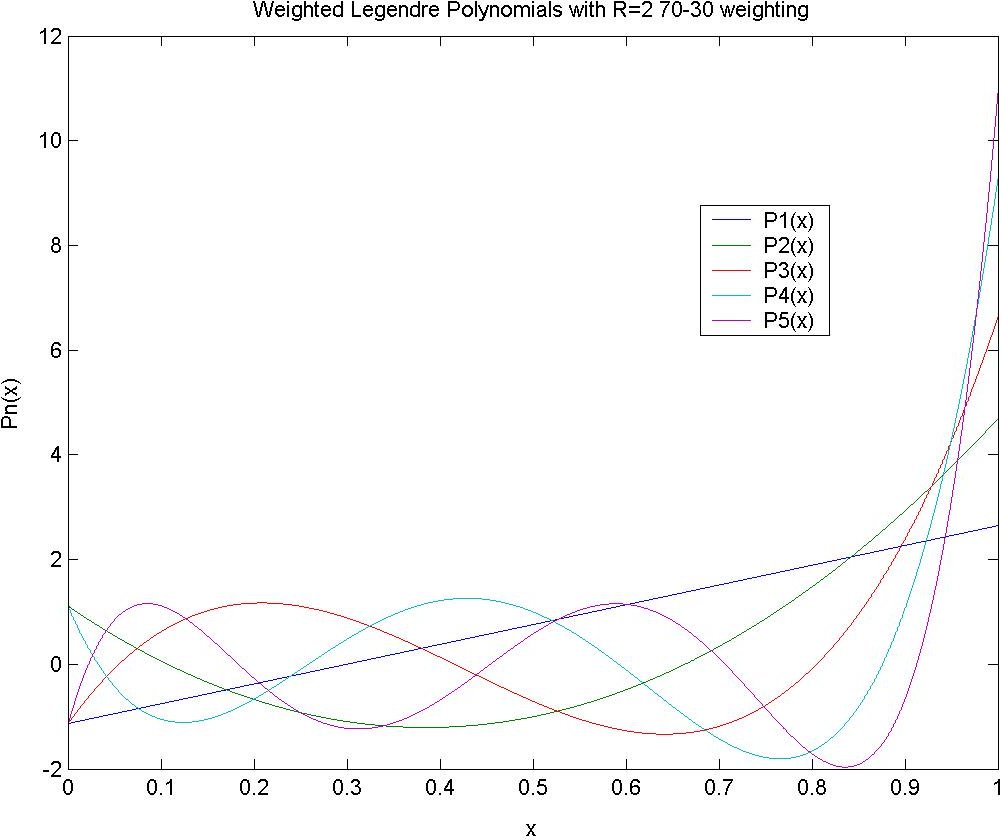}
\caption{Plot of $P_{n}(x)$  on the unit interval with $p_{1}=.7, p_{2}=.3$.}
\label{figthree3}
\end{figure}

\begin{figure}[htbp!]
\includegraphics[width=\onewidth\textwidth]{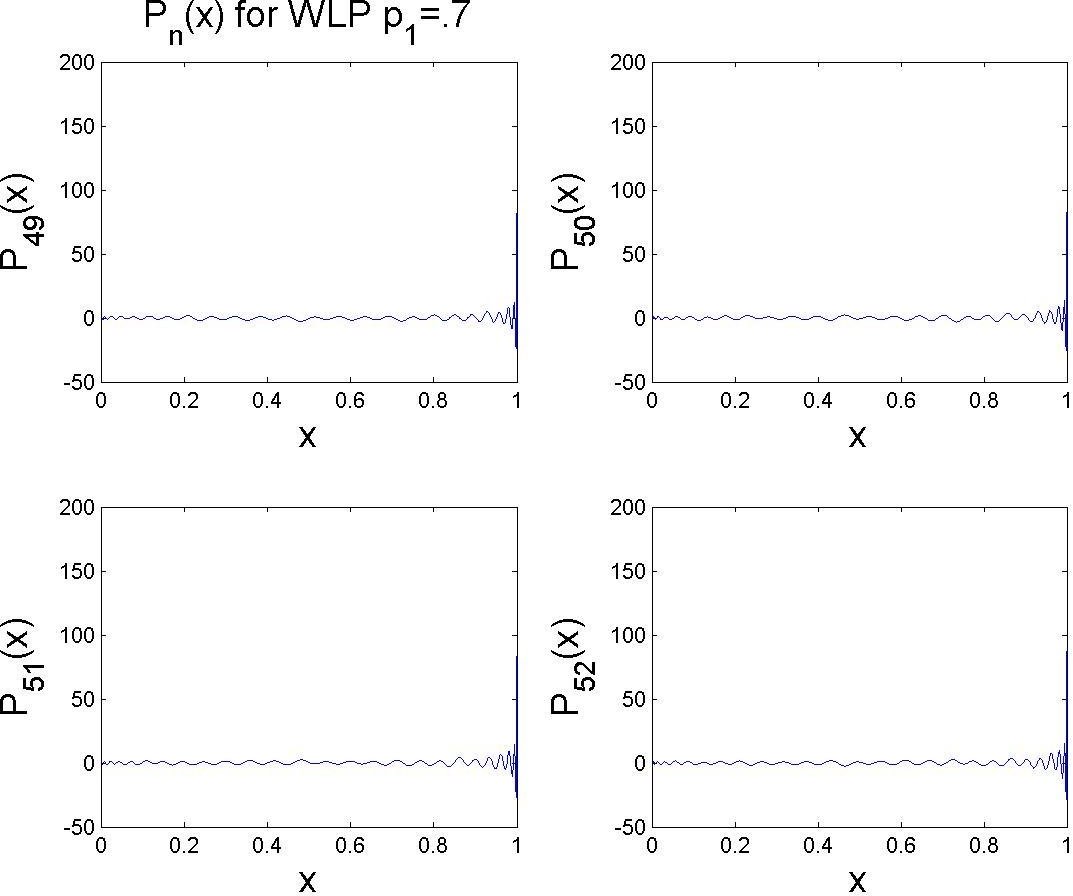}
\caption{Plot of $P_{n}(x)$  on the unit interval with $p_{1}=.7, p_{2}=.3$ for $49\leq n\leq52$.}
\label{figthree4}
\end{figure}

\begin{figure}[htbp!]
\includegraphics[width=\onewidth\textwidth]{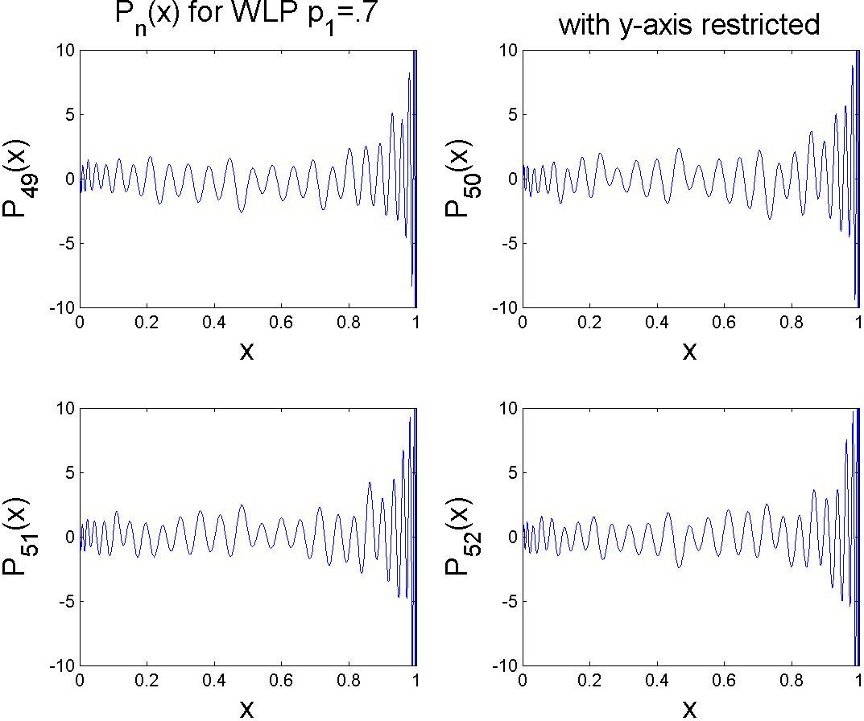}
\caption{Plot of $P_{n}(x)$  on the unit interval with $p_{1}=.7, p_{2}=.3$ for $49\leq n\leq52$ with the $y$-axis restricted to the range $-10\leq y\leq 10$.}
\label{figthree5}
\end{figure}

In Fig. \ref{figthree1} we show the graphs of $P_{n}(x)$ for $1\leq n\leq 5$ for the WLP with $p_{1}=1/2$, so these are just (rescaled versions of) the classical Legendre Polynomials, computed using G. Mantica's algorithm.  In Figs. \ref{figthree2} and \ref{figthree3} we show the same functions for $p_{1}=.6$ and $p_{1}=.7$.  Already we observe that the symmetry is broken in a decisive fashion, as these functions are much larger near $x=1$ (where the measure ``has less weight'') than near $x=0$ (where the measure ``has more weight'').  This observation is expected since the Gram-Schmidt process with respect to $\mu$ almost immediately implies the following property: $P_{n}(x)$ is the unique $n^{th}$ degree polynomial, with highest degree coefficient $d_{n}$, which minimizes the $L^{2}([0,1],\mu)$ norm.  (In the case $p_{1}=1/2$, we have $d_{n}=\prod_{i=1}^{n}\frac{1}{b_{i}}$ and $1/b_{i}=\sqrt{16-4/i^{2}}$.)  Therefore, $P_{n}$ is expected to be the smallest where $\mu$ ``has more support.''  In Fig \ref{figthree4} we show the graphs of $P_{n}(x)$ for $49\leq n\leq 52$ for $p_{1}=.7$.  The same graphs are shown in Fig. \ref{figthree5} with the $y$-axis truncated to more clearly display the structure of the functions.

\underline{Conjecture 3.1}: For fixed $p_{1}>1/2$ and any given $\epsilon\geq0$ the WLP are uniformly bounded on the interval $[0,1-\epsilon]$.

\begin{figure}[htbp!]
\includegraphics[width=\onewidth\textwidth]{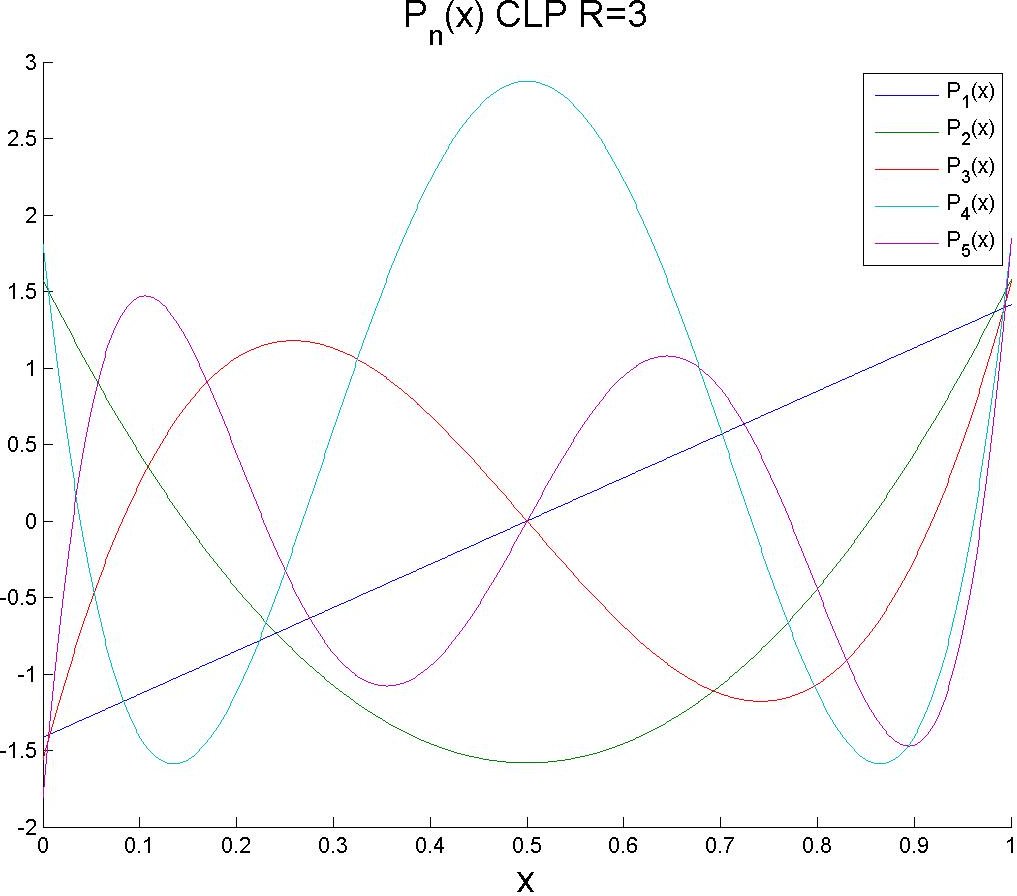}
\caption{Plot of first five Cantor Legendre Polynomials, $R=3$.}
\label{figthree6}
\end{figure}

\begin{figure}[htbp!]
\includegraphics[width=\onewidth\textwidth]{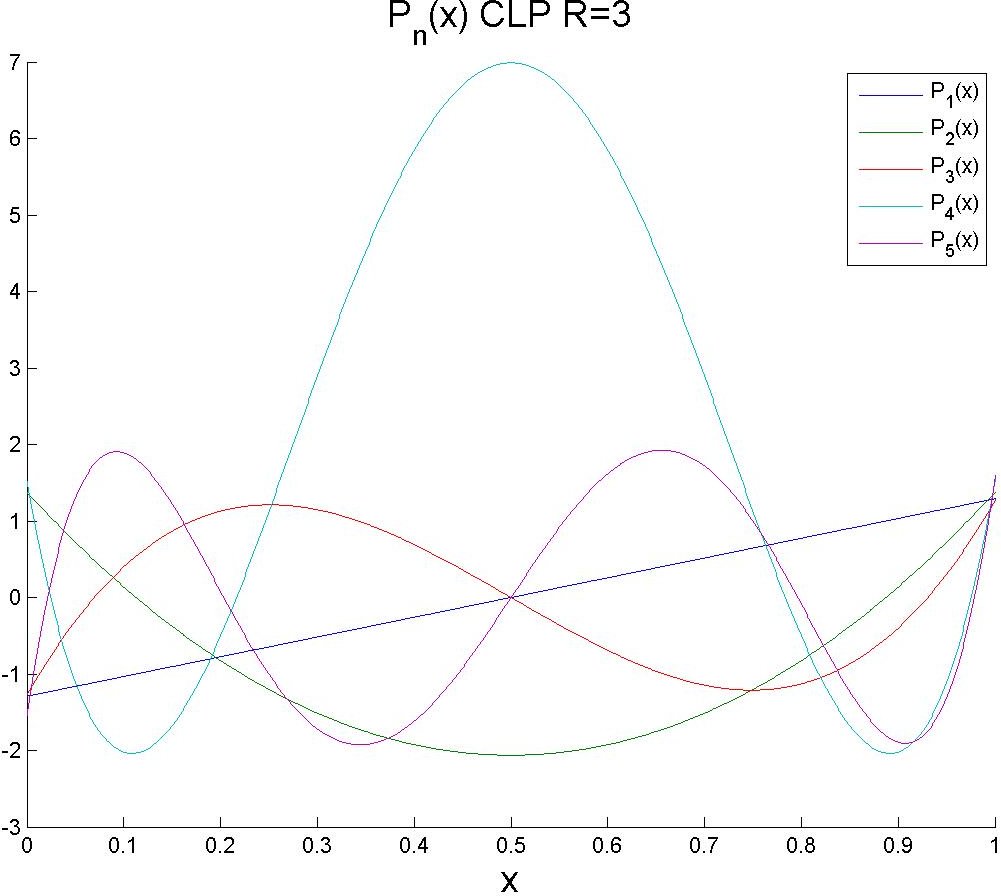}
\caption{Plot of first five Cantor Legendre Polynomials, $R=4$.  As in the Cantor Legendre Case, the polynomials have lowest absolute value where the measure has the most weight (as expected).}
\label{figthree7}
\end{figure}

\begin{figure}[htbp!]
\includegraphics[width=\onewidth\textwidth]{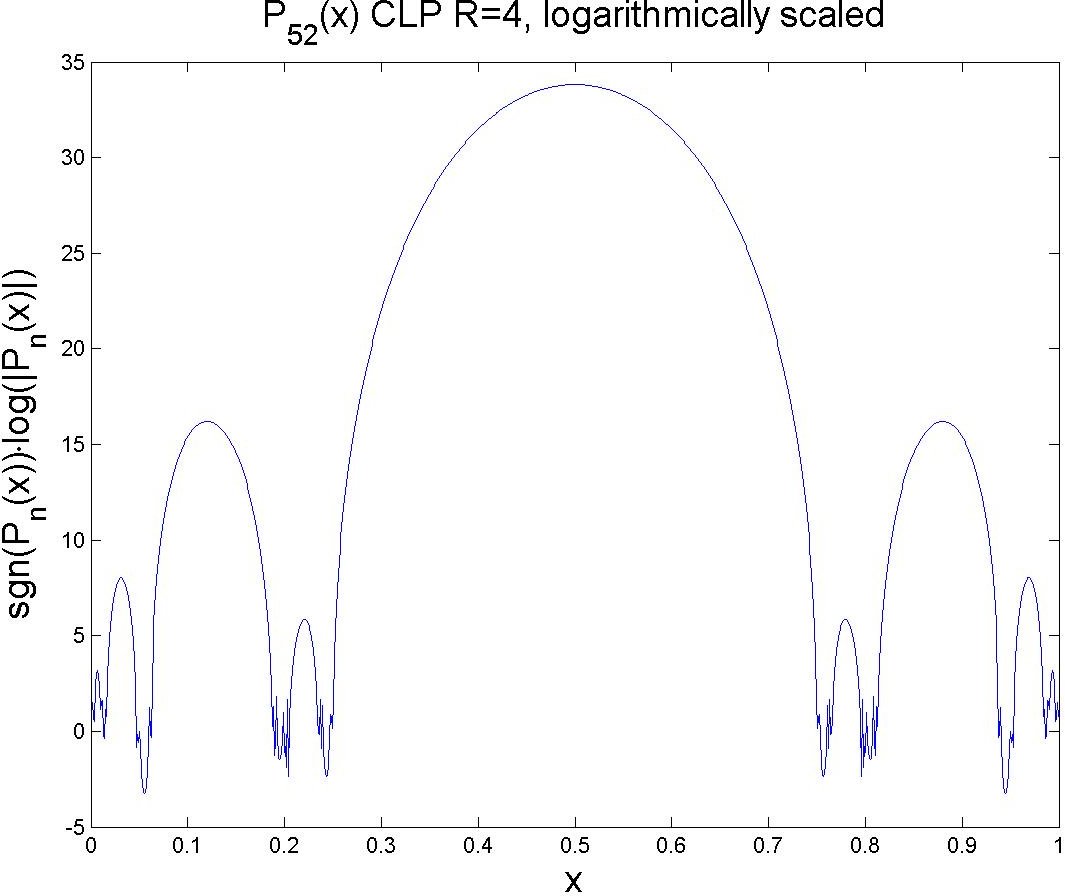}
\caption{Plot of $P_{52}(x)$ for CLP $R=4$, logarithmically scaled on the entire interval.  Notice the large values obtained by the polynomials off the Cantor set $(10^{28})$ compared to the behavior on the Cantor set (for $1\leq j\leq n$ large, $\abs{P_{n}}<20$).  The overall shape in the gap is Gaussian, and this occurs for all $n=4k$, $k\in\Z^{+}$.  For $n=4k+2$, the behavior in the gap is roughly a negative Gaussian.}
\label{figthree8}
\end{figure}

\begin{figure}[htbp!]
\includegraphics[width=\onewidth\textwidth]{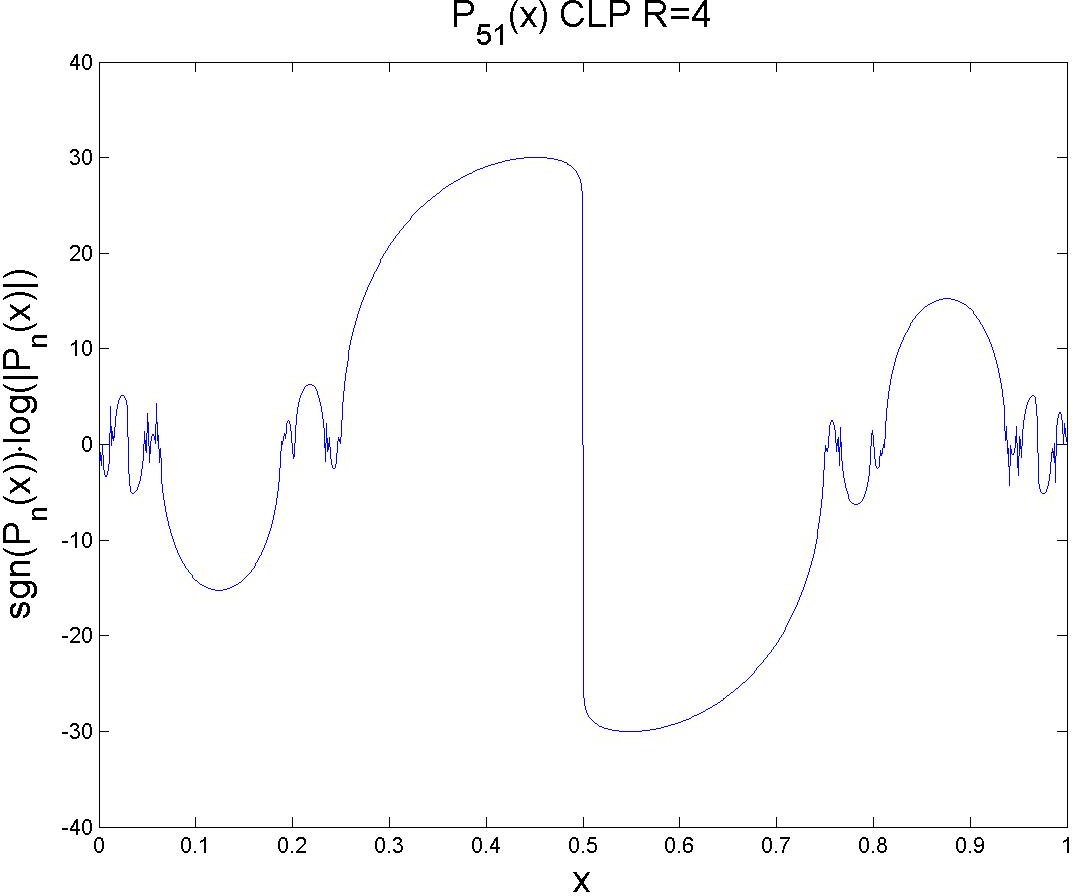}
\caption{Plot of $P_{51}(x)$ for CLP $R=4$, logarithmically scaled on the entire interval.  This is roughly the derivative of a Gaussian on the center gap.  For all $n=4k+3$ we observe this shape, and for $n=4k+1$, the negative of this shape.}
\label{figthree9}
\end{figure}

\begin{figure}[htbp!]
\includegraphics[width=\onewidth\textwidth]{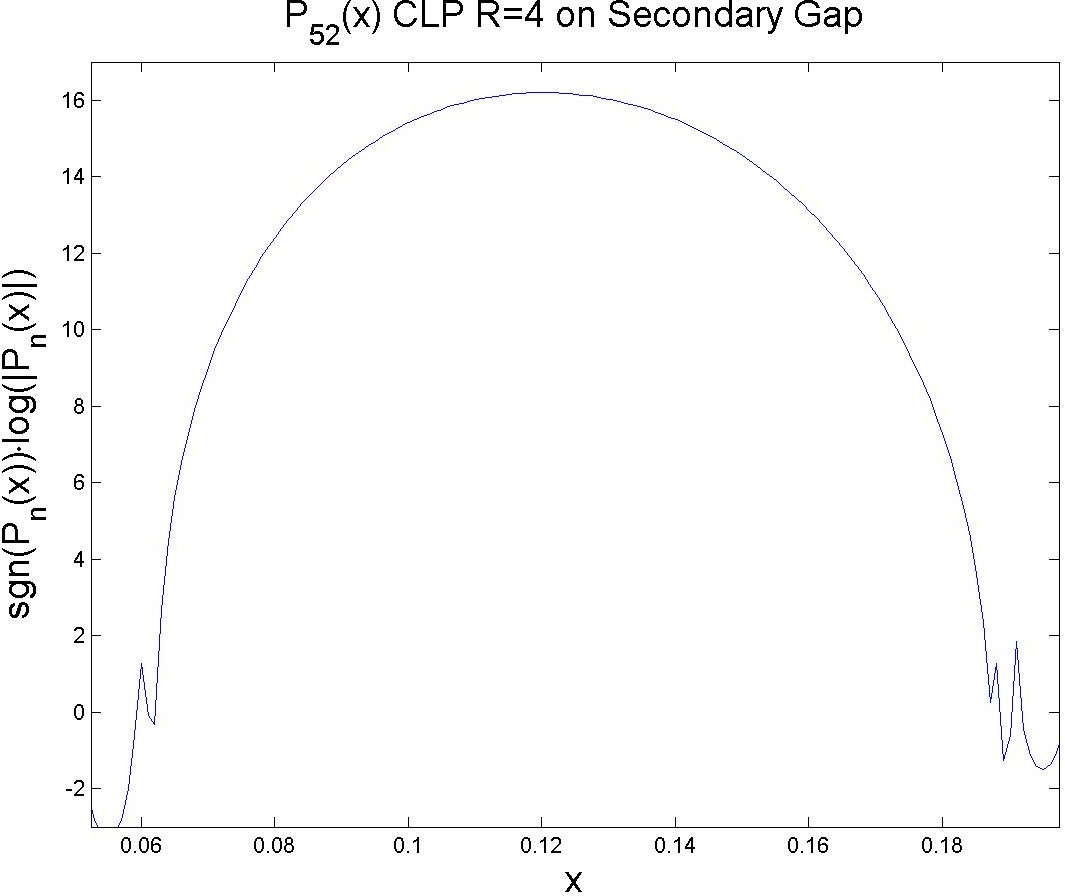}
\caption{Plot of $P_{52}(x)$ for CLP $R=4$, logarithmically scaled on the secondary gap $[1/R,(R-1)/R]=[1/16,3/16]$.  Again, we observe a Gaussian shape.}
\label{figthree10}
\end{figure}

Next we look at the CLP case, where certain features of the $P_{n}(x)$ increase in complexity.  In Figs. \ref{figthree6} and \ref{figthree7} we show the graphs of $P_{n}(x)$ for $1\leq n\leq 5$ on the whole interval for two different choices of $R$.  In Fig. \ref{figthree8} we show the graph of $P_{52}(x)$ on the whole interval for $R=4$.  Note that the values on the central gap are so large that no information about the graph on the complement of the central gap is discernable using linear scaling in the $y$-axis.  This is typical of $P_{n}(x)$ for large even values of $n=2k$ and any $R>2$ (though the sign of $P_{n}(x)$ changes sign between $n=4j+2$ and $n=4j$).  To view these polynomials more effectively, we scale the $y$-axis logarithmically, but also multiply by the sign of $P_{n}(x)$.  That is, we plot $\sign(P_{n}(x))\log(\abs{P_{n}(x)})$ vs $x$.  In Fig. \ref{figthree9} we show the graph of $P_{51}(x)$ for $R=4$, which is typical of $P_{n}(x)$ for $n$ odd and large.  With the same logarithmic scaling, Fig. \ref{figthree10} shows $P_{52}(x)$ for $R=4$ restricted to the gap $[1/R^{2},]=[1/16,3/16]$.  This behavior is typical of $P_{n}(x)$ for large $n$.  Note that the maximum value is quite large, but still it is small relative to the maximum value on the central gap.  We will discuss this behavior more in Section \ref{secclp}.

\begin{figure}[htbp!]
\includegraphics[width=\onewidth\textwidth]{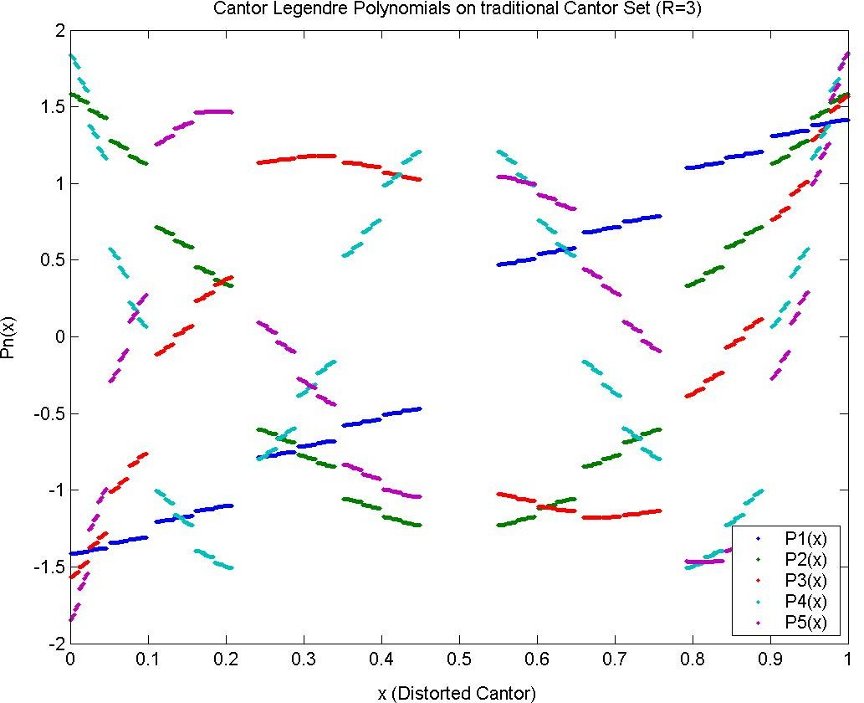}
\caption{Plot of first five Cantor Legendre Polynomials, $R=3$.}
\label{figthree11}
\end{figure}

\begin{figure}[htbp!]
\includegraphics[width=\onewidth\textwidth]{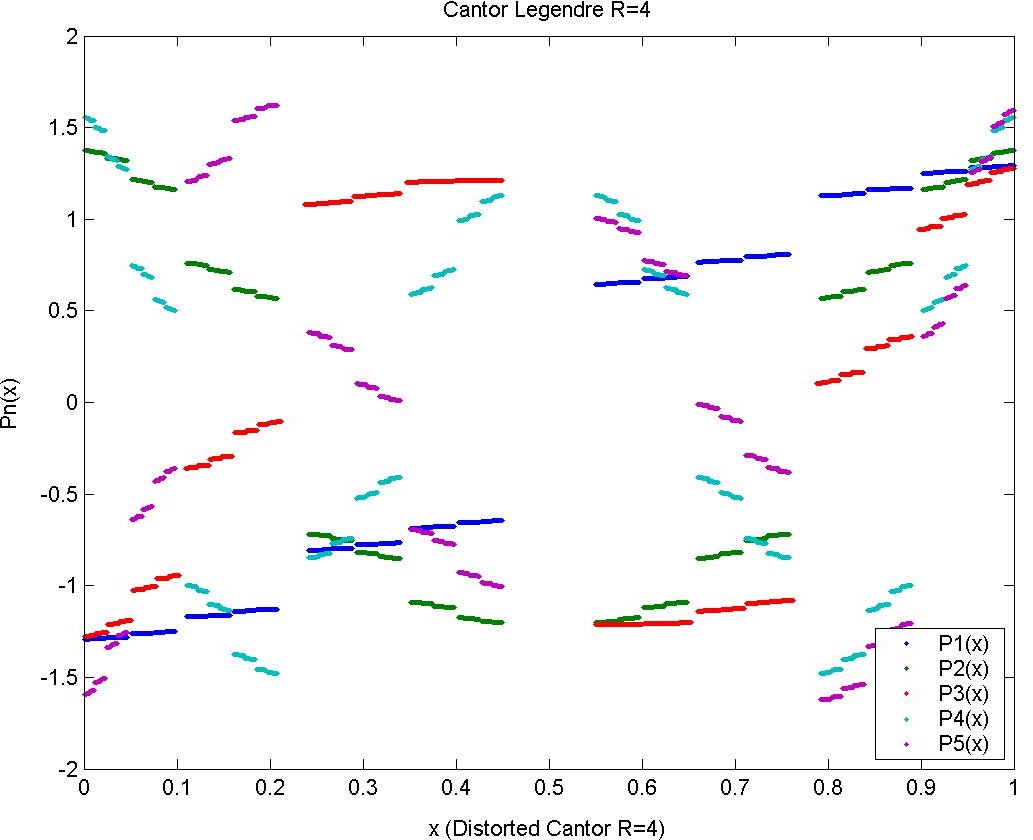}
\caption{Plot of first five Cantor Legendre Polynomials, $R=4$}
\label{figthree12}
\end{figure}

\begin{figure}[htbp!]
\includegraphics[width=\onewidth\textwidth]{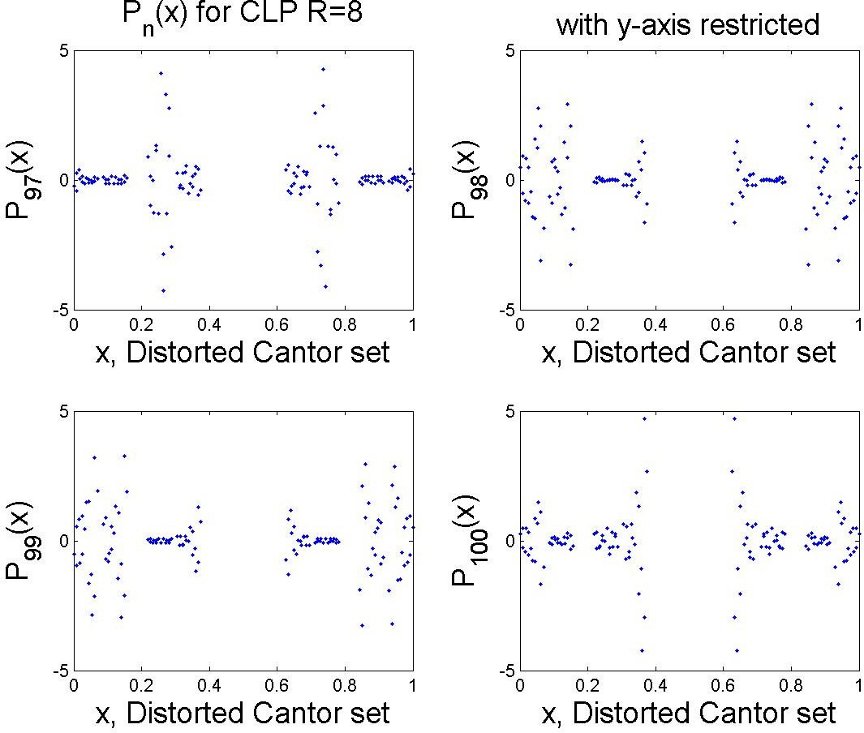}
\caption{Plot of some larger $n$ Cantor Legendre Polynomials, $R=8$, $97\leq n\leq100$.}
\label{figthree13}
\end{figure}

Next, we graph CLP on the Cantor set $C_{R}$ using a distorted Cantor set $\widetilde{C}_{R,\epsilon}$ for the $x$-axis.  In Figs. \ref{figthree11} and \ref{figthree12} we display the same functions as in Figs. \ref{figthree6} and \ref{figthree7}.  Already we see that the values of $P_{n}(x)$ are considerably smaller on the Cantor set.  In Fig. \ref{figthree13} we show $P_{n}(x)$ restricted to the Cantor set for $97\leq n\leq 100$ and $R=8$.  We will comment in detail about some of the structure of those graphs in Section \ref{secapprox}.  Many more examples may be viewed at \cite{owrutsky05}.

\clpg
\section{Dirichlet Kernels}
\label{secdir}

For a general function $f$ in $L^{2}(d\mu)$, we can expand it as a series
\begin{equation}\label{four15}
\sum_{n=0}^{\infty}a_{n}P_{n}(x),
\end{equation}
where the coefficients are given by
\begin{equation}\label{four16}
a_{n}=\int f(y)P_{n}(y)\,d\mu(y).
\end{equation}
The partial sums
\begin{equation}\label{four17}
\sum_{n=0}^{N}a_{n}P_{n}(x)
\end{equation}
may be represented as an integral
\begin{equation}\label{four18}
\int D_{N}(x,y)f(y)\,d\mu(y)
\end{equation}
where the Dirichlet kernel is given by
\begin{equation}\label{four19}
D_{n}(x,y)=\sum_{n=0}^{N}P_{n}(x)P_{n}(y)
\end{equation}
The partial sums in Eq. (\ref{four17}) converge to $f$ in $L^{2}$ norm, but to get better convergence we need to know more about the Dirichlet kernel.  In view of related results in \cite{strichartz932,strichartz94}, one might hope that there exists a sequence $\{N_{k}\}$ along which the partial sums converge uniformly if $f$ is continuous.  This would follow by standard approximate identity arguments if we could show that the quantity
\begin{equation}\label{four20}
\int\abs{D_{N_{k}}(x,y)}\,d\mu(y)
\end{equation}
is uniformly bounded, and that
\begin{equation}\label{four21}
\lim_{k\to\infty}\int_{\abs{y-x}\geq\epsilon}\abs{D_{N_{k}}(x,y)}\,d\mu(y)=0
\end{equation}
for all $\epsilon>0$.  While we have no insight on how to establish Eq. (\ref{four20}), we can say something about Eq. (\ref{four21}), thanks to the Christoffel-Darboux formula
\begin{equation}\label{four22}
D_{N}(x,y)=\frac{r_{N+1}}{x-y}(P_{N+1}(x)P_{N}(y)-P_{N}(x)P_{N+1}(y))\mbox{ , for }x\neq y
\end{equation}
(It is easy to derive Eq. (\ref{four22}) by multiplying Eq. (\ref{four19}) by $(x-y)$ and using the $3$-term recursion relation Eq. (\ref{one3}).)  Assuming that the polynomials $P_{n}(x)$ are uniformly bounded on the support of $\mu$, if we could find a sequence $N_{k}$ such that $r_{N_{k+1}}\to0$, then Eq. (\ref{four21}) follows from Eq. (\ref{four22}).

It appears from our data that for CLP polynomials there exist indices $n$ such that $r_{n}$ is close to zero, but there is no evidence for a sequence tending to zero.  This means that there will be some Dirichlet kernels that seem very concentrated near the diagonal (those with $r_{N+1}$ close to zero), but it is unlikely that we can improve this behavior indefinitely.

\begin{figure}[htbp!]
\includegraphics[width=\onewidth\textwidth]{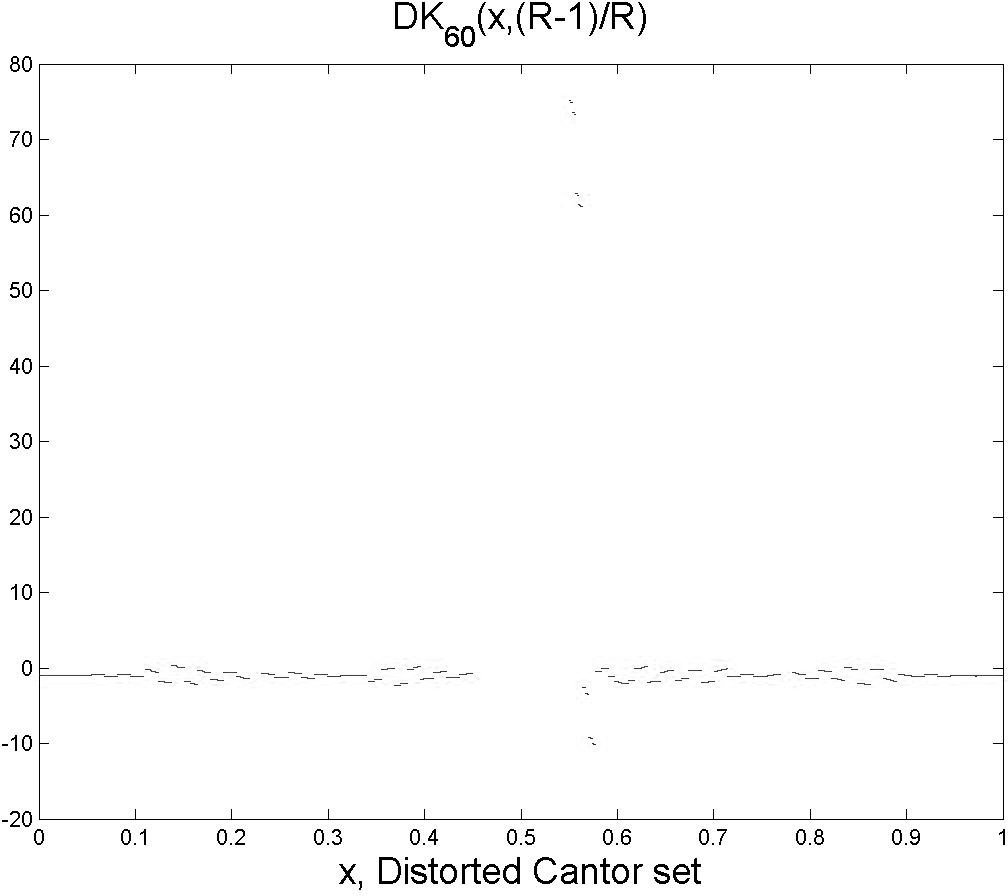}
\caption{Dirichlet kernel for CLP with $n=60$ and $R=8$ centered at the left most point in the right half of the Cantor set.  We see that the graph is reasonably small away from the center.}
\label{figfour1}
\end{figure}

\begin{figure}[htbp!]
\includegraphics[width=\onewidth\textwidth]{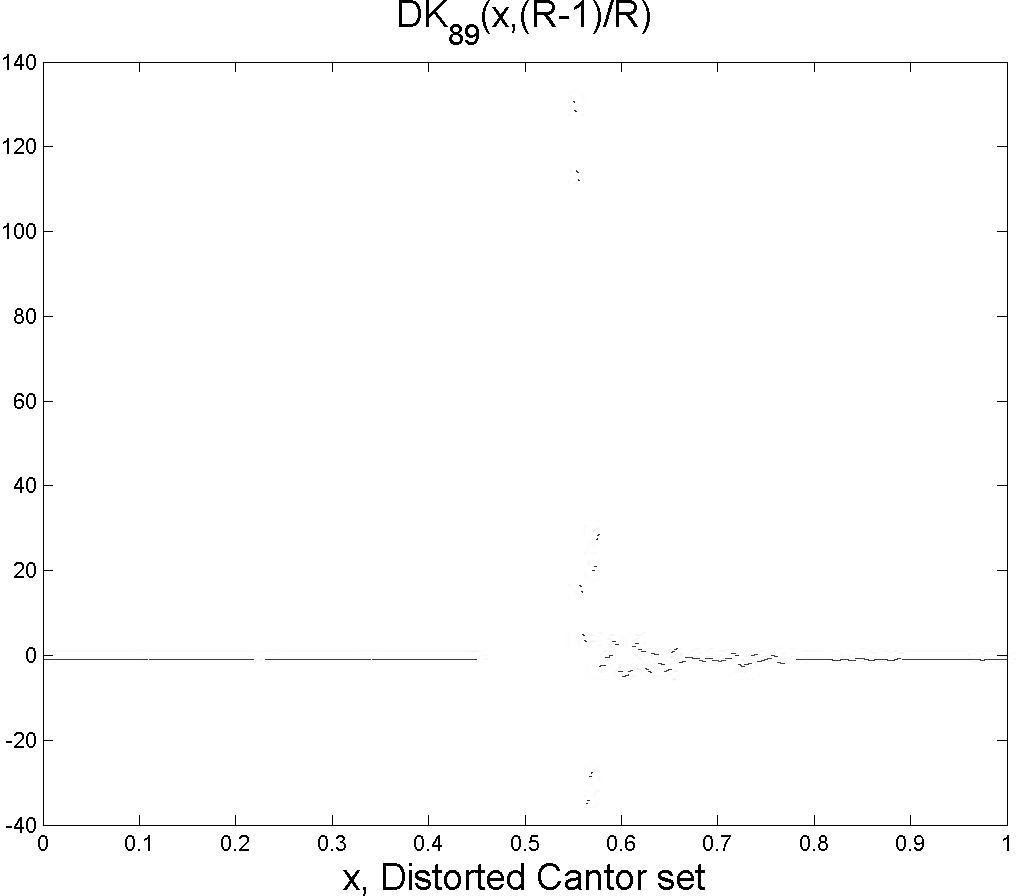}
\caption{Dirichlet kernel for CLP with $n=89$ and $R=8$ centered at the left most point in the right half of the Cantor set.  Again, the kernel is small away from its center.}
\label{figfour2}
\end{figure}

\begin{figure}[htbp!]
\includegraphics[width=\onewidth\textwidth]{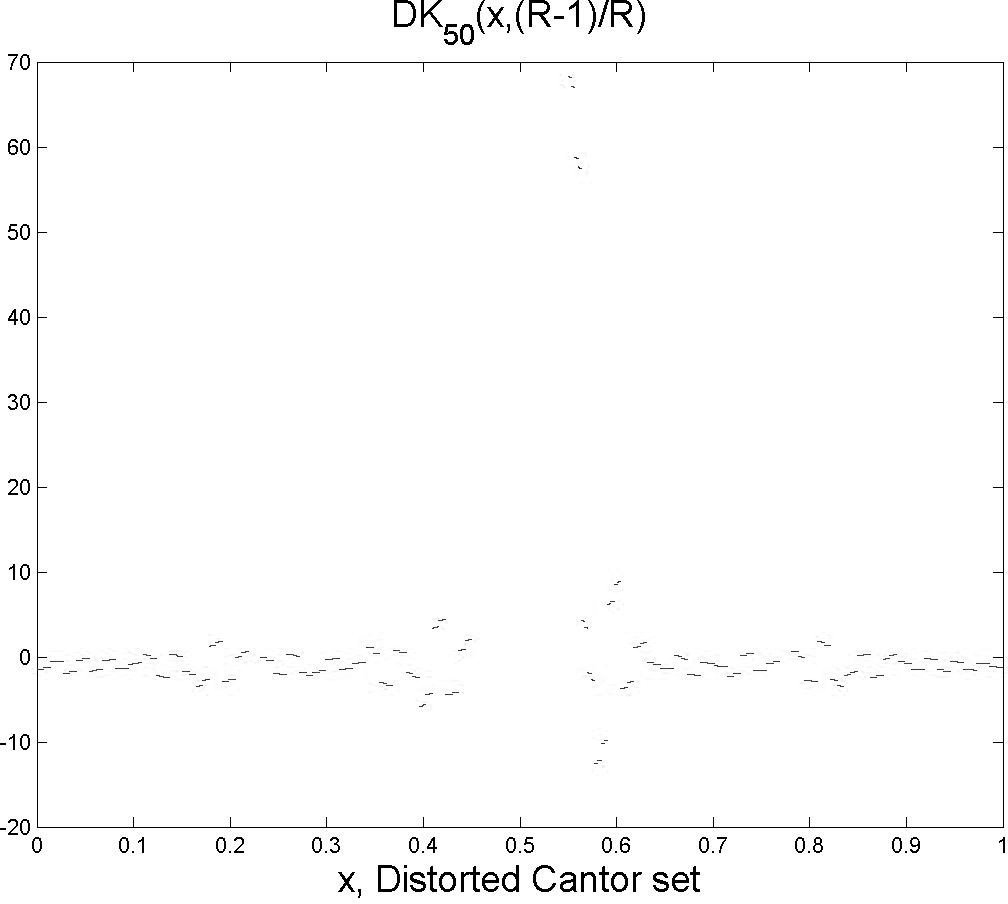}
\caption{Dirichlet kernel for CLP with $n=50$ and $R=8$ centered at the left most point in the right half of the Cantor set.  This kernel takes much larger values away from the center than the $n=60$ case in Fig. \ref{figfour1}.}
\label{figfour3}
\end{figure}

\begin{figure}[htbp!]
\includegraphics[width=\onewidth\textwidth]{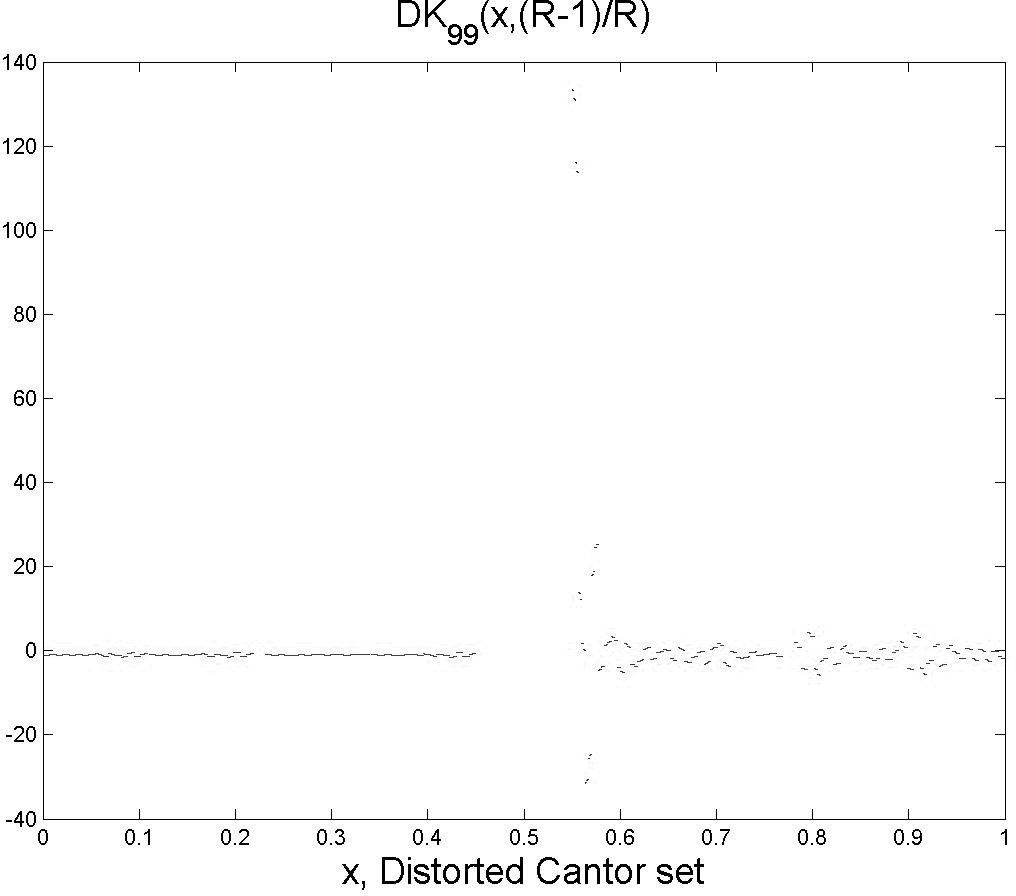}
\caption{Dirichlet kernel for CLP with $n=99$ and $R=8$ centered at the left most point in the right half of the Cantor set.}
\label{figfour4}
\end{figure}

Figs. \ref{figfour1} and \ref{figfour2} illustrate Dirichlet kernels $D_{N}(\cdot,y)$ for fixed $y$  that are not concentrated near $x=y$, while Figs. \ref{figfour3} and \ref{figfour4} illustrate Dirichlet kernels that are moderately well concentrated.

\clpg
\section{Approximate Equalities for CLP Polynomials}
\label{secapprox}

\begin{figure}[htbp!]
\includegraphics[width=\onewidth\textwidth]{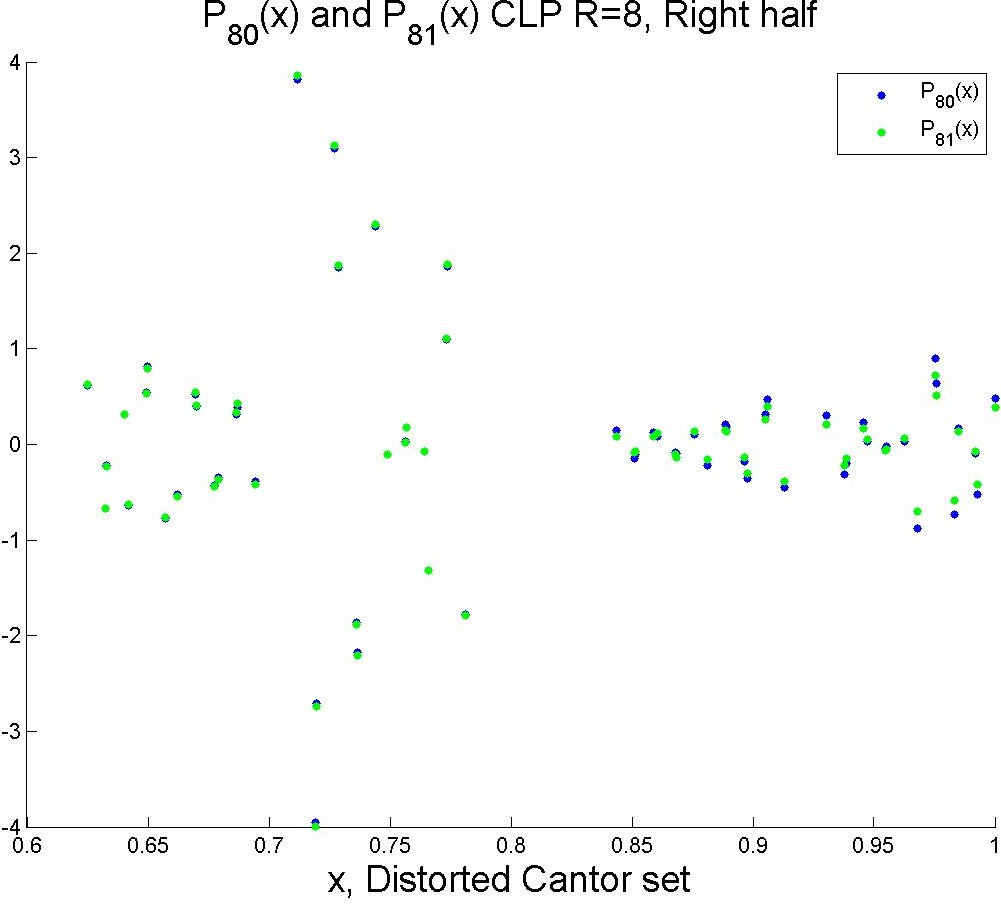}
\caption{$P_{80}(x)$ and $P_{81}(x)$, CLP $R=8$}
\label{figfive1}
\end{figure}

\begin{figure}[htbp!]
\includegraphics[width=\onewidth\textwidth]{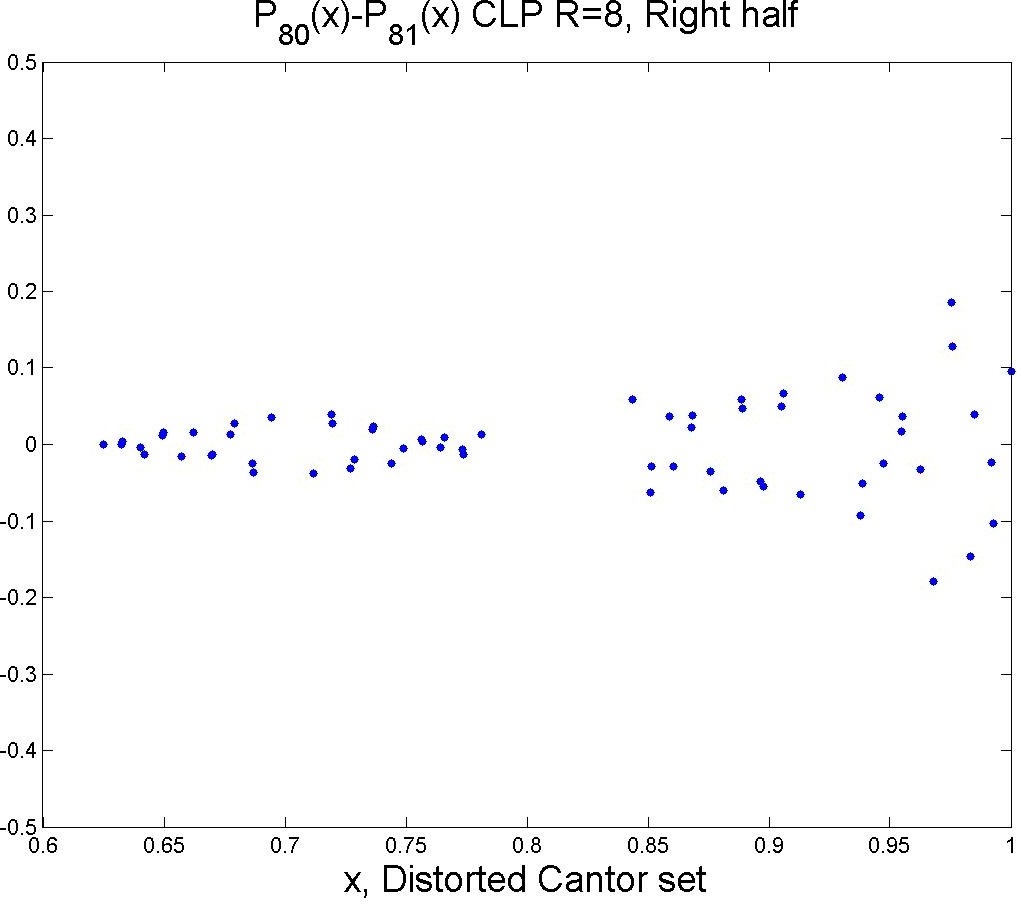}
\caption{$P_{80}(x)-P_{81}(x)$, CLP $R=8$}
\label{figfive2}
\end{figure}

\begin{figure}[htbp!]
\includegraphics[width=\onewidth\textwidth]{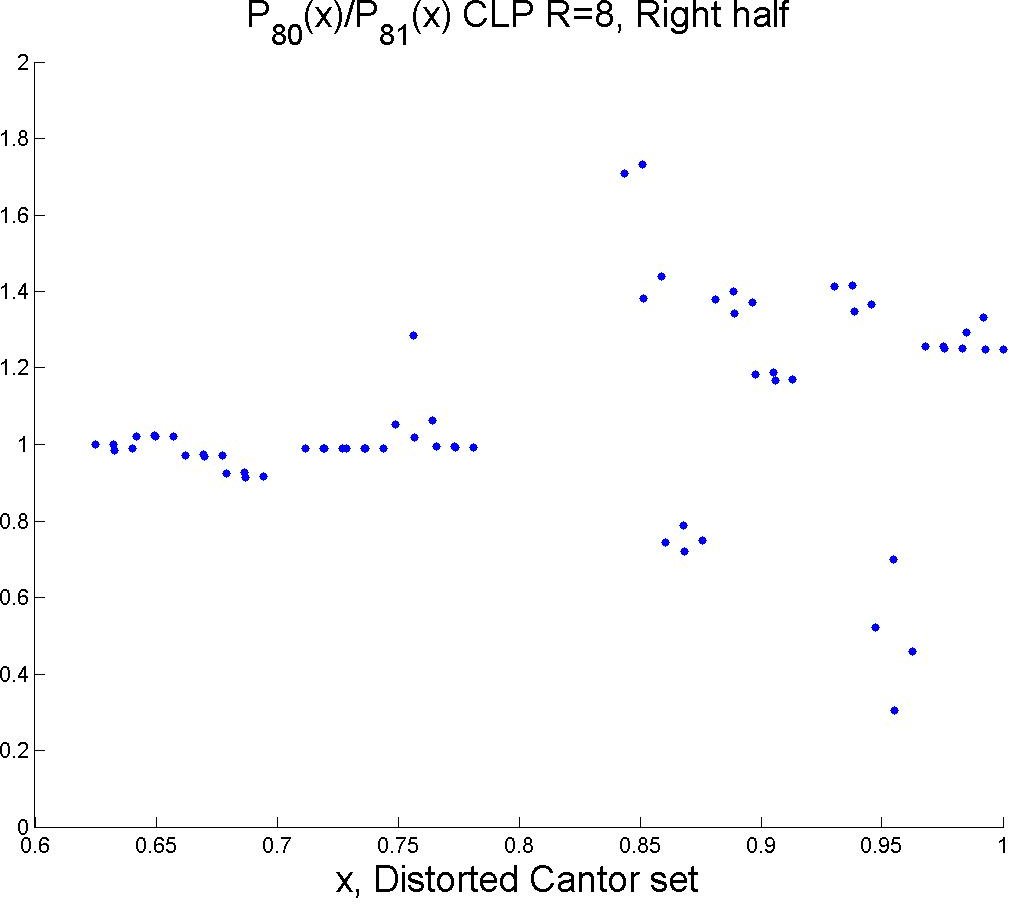}
\caption{$P_{80}(x)/P_{81}(x)$, CLP $R=8$}
\label{figfive3}
\end{figure}

In this section we discuss two types of approximate equalities in the CLP case. First we note that by symmetry, $P_{n}(x)$ is even under the reflection $x\mapsto1-x$ when $n$ is even and odd when $n$ is odd.  Nevertheless, the plots of $P_{2n}(x)$ and $P_{2n+1}(x)$ appear very similar (see $n=98,99$ in Fig. \ref{figthree13}).  This is especially true on the right half of the Cantor set.  In Fig. \ref{figfive1} we graph $P_{50}(x)$ and $P_{51}(x)$ for $R=8$.  The qualitative similarity is striking, but it is difficult to quantify.  Fig. \ref{figfive2} shows the graph of the difference, and Fig. \ref{figfive3} shows the graph of the ratio.  We can give a rough explanation using Eq. (\ref{one4}), which for even values may be written
\begin{equation}\label{five1}
P_{2n+1}(x)=\frac{x-1/2}{r_{2n+1}}P_{2n(x)}-\frac{r_{2n}}{r_{2n+1}}P_{2n-1}(x).
\end{equation}
For large values of $R$, we may have $r_{2n}$ close to zero and $r_{2n+1}$ close to $1/2$.  So the second term on the right side of Eq. (\ref{five1}) is close to zero.  But on the right half of the Cantor set, $x$ is close to one, so the coefficient $(x-1/2)/r_{2n+1}$ is close to one.  Therefore, Eq. (\ref{five1}) says that $P_{2n+1}(x)\approx P_{2n}(x)$ on the right half of $C_{R}$.  By the odd-even behavior, we have $P_{2n+1}(x)\approx -P_{2n}(x)$ on the left half of $C_{R}$.  Since we see a qualitative similarity of the plots on all of $C_{R}$, we must attribute these approximate equalities to an approximate reflectional symmetry across the $x$-axis of the graphs of all $P_{n}(x)$.  This is roughly apparent in Fig. \ref{figthree13}, but it does not hold up to close inspection.  In particular, the graphs of $P_{2n+1}(x)$ and $P_{2n+1}(1-x)$ are not that close (of course $P_{2n+1}(1-x)=-P_{2n+1}(x)$ exactly).

\begin{figure}[htbp!]
\includegraphics[width=\onewidth\textwidth]{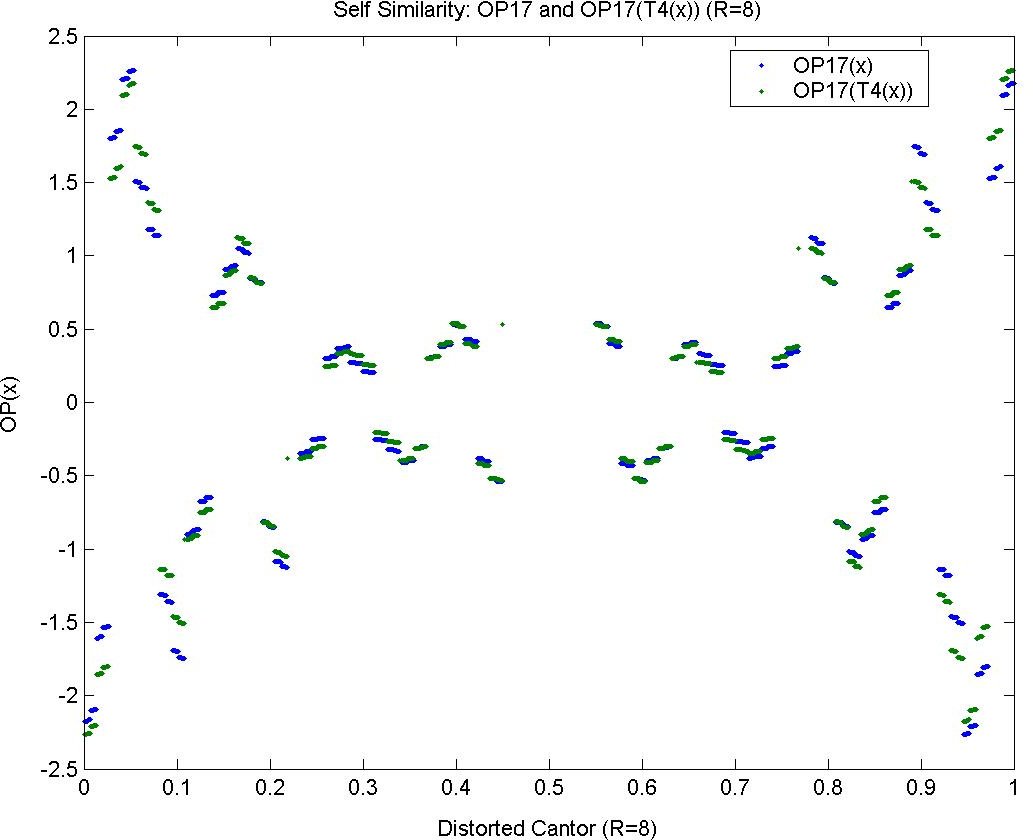}
\caption{$P_{17}(x)$ and $P_{17}(T_{3}(x))$, $R=8$}
\label{figfive4}
\end{figure}

The second type of approximate equality refines the above idea.  If we express the points of $C_{R}$ as infinite binary decimals, then $x\mapsto (1-x)$ simply interchanges all digits in the binary expansion.  Let $T_{m}(x)$ denote the map $C_{R}\to C_{R}$ that interchanges the first $m$ binary digits, leaving all other digits unchanged.  In other words, $T_{m}$ permutes the $2^{m}$ Cantor subsets of level $m$ by reversing the order of the subsets.  In Fig. \ref{figfive4} we show the plots of $P_{17}(x)$ and $P_{17}(T_{3}(x))$ for $R=8$.  There is clearly a strong qualitative fit, but the agreement of the two functions is not very close numerically.  The same pattern persists for $P_{2^{m}+1}(x)$ and $P_{2^{m}+1}(T_{m}(x))$ for all values of $R$.  At present we have no explanation for this phenomenon.

\clpg
\section{CLP Polynomials on Gaps}
\label{secclp}

\begin{figure}[htbp!]
\includegraphics[width=\onewidth\textwidth]{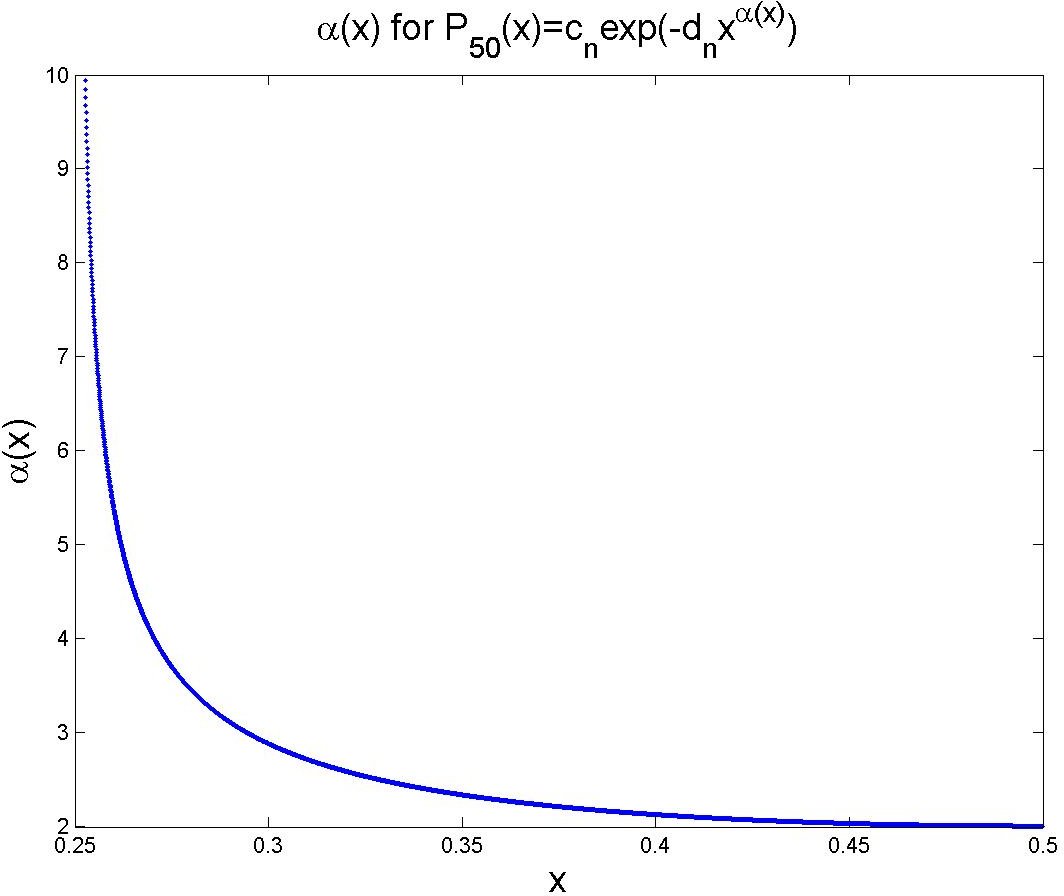}
\caption{$\alpha(x)$ for CLP $R=4$, $n=50$}
\label{figsix1}
\end{figure}

We saw in Figs. \ref{figthree8} and \ref{figthree10} that the graph of $P_{2n}(x)$ on the gaps in $C_{R}$ is approximately Gaussian, for $n$ large.  In reality, we find that $P_{2n}|_{[1/R,(R-1)/R]}(x)=c_{2n}\exp{(-d_{2n}x^{\alpha(x)})}$ where $\alpha(x)=\alpha_{2n}(x)=2$ at $x=1/2$, and $\alpha(x)\approx2$ for $x$ around $1/2$.  Fig. \ref{figsix1} shows $\alpha(x)=\alpha_{2n}(x)$ for $n=50$, CLP $R=4$.  Due to symmetry, we only view $\alpha$ on the interval $[1/R,1/2]$.  From the data, we conjecture that $\alpha_{2n}$ converges as $n\to\infty$ on any closed interval of the form $[1/R+\epsilon,1/2]$ to a function resembling that in Fig. \ref{figsix1}.

The data also indicate exponential growth (in $n$) of $P_{2n}(x)$ on the central gap $[1/R,(R-1)/R]$.  This exponential growth results from the behavior of Eq. (\ref{one4}) at $x=1/2$ (with $A_{n}=1/2$ for all $n$) as follows:
\begin{equation}\label{six1}
P_{2n}(1/2)=(-1)^{n}\prod_{j=1}^{n}(r_{2j-1}/r_{2j}).
\end{equation}
Numerically, we find that $(r_{2j-1}/r_{2j})$ is bounded above $1$ with few exceptions.  In the CLP $R=2.5$ case, there are $134$ instances that $(r_{2j-1}/r_{2j})\leq1$ for $j=1,\ldots,5000$, in the $R=3$ case there is one, and for the $R=4$ case there are none.   Thus, as $R$ increases, $(r_{2j-1}/r_{2j})$ is bounded above $1$ more consistently, as we have seen in the $R=8$ case in Figs. \ref{figtwo6} and \ref{figtwo7}.  Eq. (\ref{six1}) therefore (generally) gives exponential growth in $n$ for $\abs{P_{2n}(1/2)}$.

We then summarize the behavior of $P_{2n}(x)$ on the gap $[1/R,(R-1)/R]$ as follows: $c_{2n}=P_{2n}(1/2)$ grows exponentially (according to Eq. \ref{six1}), $d_{2n}$ grows linearly, and $\alpha_{2n}(x)$ appears to converge (to something resembling Fig. \ref{figsix1}).  In fact, we can say more: for $k\in\Z^{+}$, write
\begin{equation}\label{six2}
P_{j}(x)=\sum_{k=0}^{j}\alpha_{j,k}(x-1/2)^{k}.
\end{equation}
Eq. \ref{six1} gives all $\alpha_{j,0}$ (observe that $\alpha_{2n,0}=c_{2n}$ and $\alpha_{2n+1,0}=0$ for all $n$), and then we can again use Eq. \ref{one4} to compute all of the $\alpha_{j,k}$.  (Here we read Eq. \ref{one4} as an equation of polynomials.)  If we view the coefficients $\{\alpha_{j,k}\}$ as a (lower triangular) matrix, then we see that the left column $\{\alpha_{j,0}\}$ and the diagonal $\{\alpha_{j,j}\}$ determine all of the other $\alpha_{j,k}$, via Eq. \ref{one4}.  Thus, the behavior of the left column and diagonal (and the $r_{j}$ coefficients) dictates the behavior of the other $\alpha_{j,k}$.

\begin{figure}[htbp!]
\includegraphics[width=\onewidthh\textwidth]{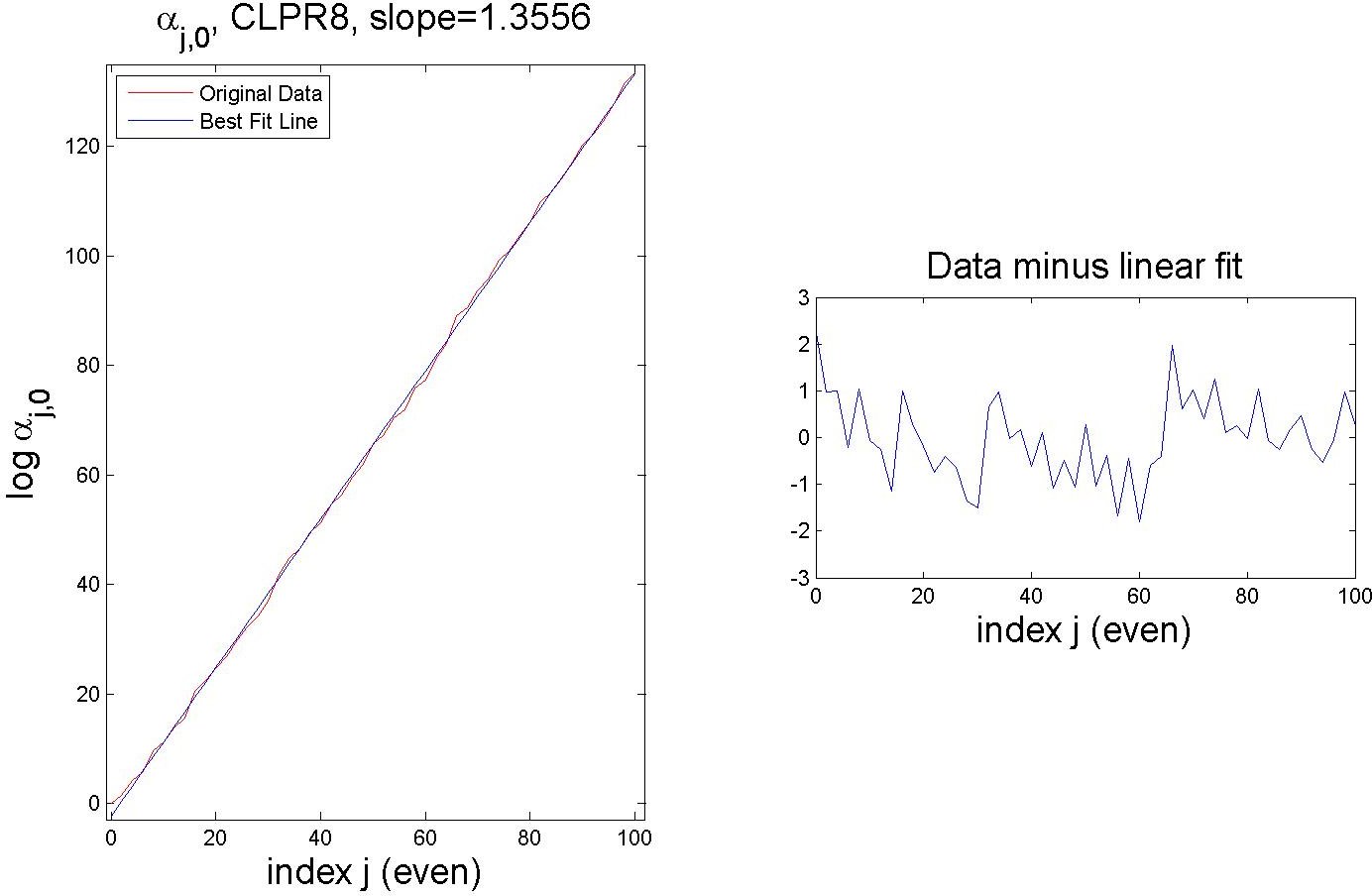}
\caption{Logarithmically scaled left column $\{\alpha_{j,0}\}$ of polynomial coefficients, with error, CLP $R=8$.  Here we only plot $j$ even since for odd $j$, $\alpha_{j,0}=0$.}
\label{figsix1.1}
\end{figure}

\begin{figure}[htbp!]
\includegraphics[width=\onewidthh\textwidth]{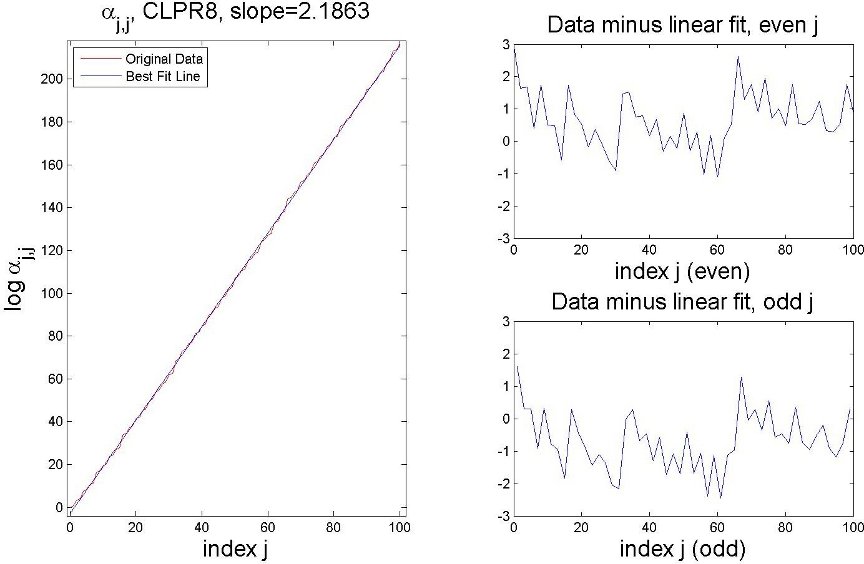}
\caption{Logarithmically scaled diagonal $\{\alpha_{j,j}\}$ of polynomial coefficients, with error, CLP $R=8$.}
\label{figsix1.2}
\end{figure}

\begin{figure}[htbp!]
\includegraphics[width=\onewidthh\textwidth]{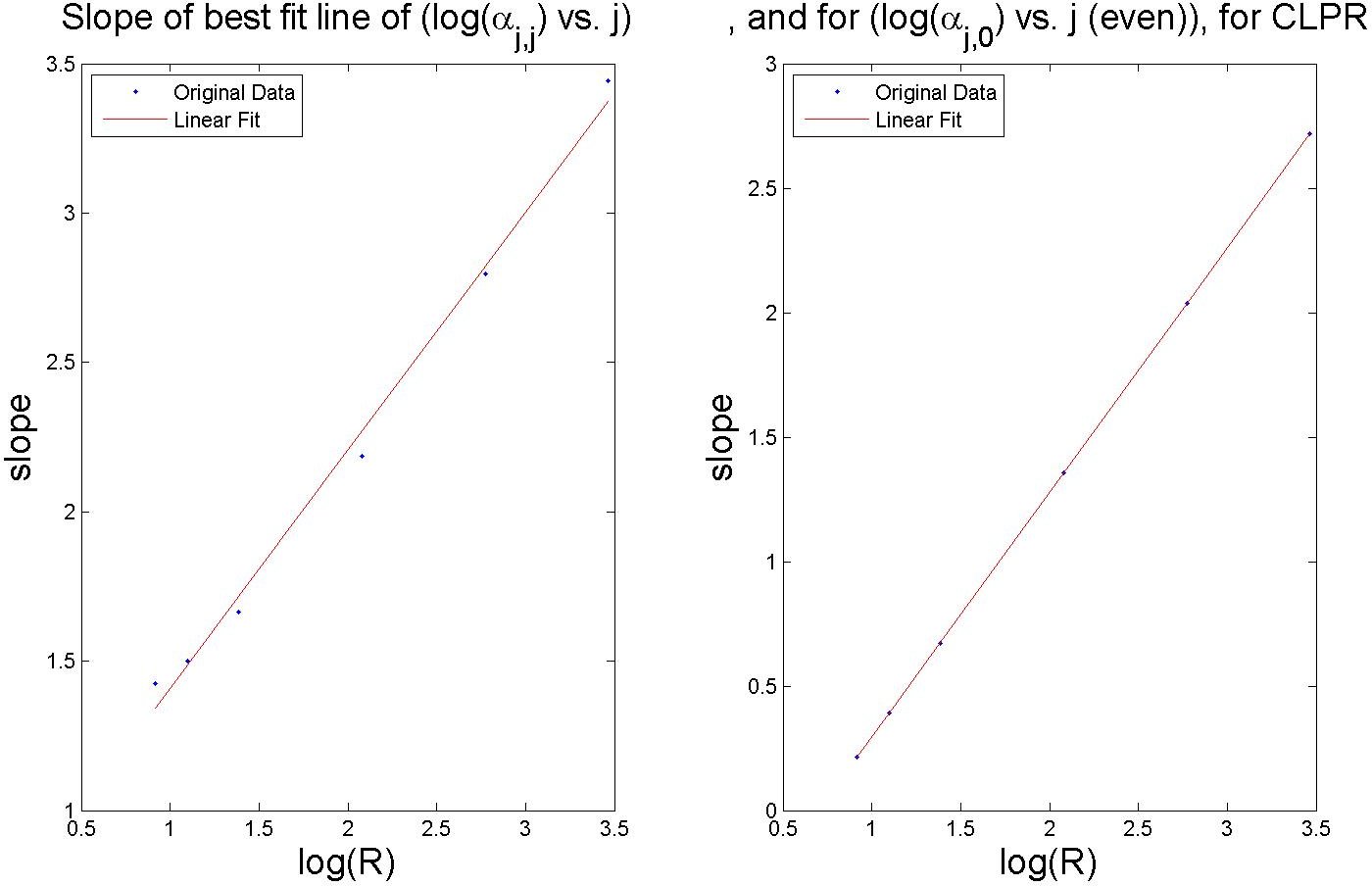}
\caption{Slopes of best fit lines for $\alpha_{j,j}$ and $\alpha_{j,0}$ (as in Figs. \ref{figsix1.1} and \ref{figsix1.2}), for $R\in\{2.5,3,4,8,16,32\}$.}
\label{figsix1.3}
\end{figure}

The left column and diagonal grow log-linearly, as we see in Figs. \ref{figsix1.1} and \ref{figsix1.2} and Table \ref{tablesix1}.  Note the similarity in the three error plots of Figs. \ref{figsix1.1} and \ref{figsix1.2}.  If we label the (top) even indexed error $\Phi_{even}$, the (bottom) odd error $\Phi_{odd}$, and the (even) error function from Fig. \ref{figsix1.1} $\Psi_{even}$, we have the following approximate equalities for $k=0,\ldots,50$: $\Phi_{even}(2k)\approx\Phi_{odd}(2k+1)+1.36\approx\Psi_{even}(2k)+.67$.  Also, from Table \ref{tablesix1} we see that $\sign(\alpha_{j,k})=-\sign(\alpha_{j,k+2})$, for appropriate $j,k$.  For the column and the diagonal of the matrix of $\alpha_{j,k}$, the slope of each best fit line varies with $\log R$, as we see in Fig. \ref{figsix1.3}.

\begin{table}[htbp!]
\resizebox{14cm}{!}{
\begin{tabular}{r|ccccccccccc}
$\alpha_{j,k}/10^4$ & $0$ & $1$ & $2$ & $3$ & $4$ & $5$ & $6$ & $7$ & $8$ & $9$ & $10$\\
\hline
0 & 0.0001 & 0 & 0 & 0 & 0 & 0 & 0 & 0 & 0 & 0 & 0 \\
1 & 0 & 0.00022678 & 0 & 0 & 0 & 0 & 0 & 0 & 0 & 0 & 0 \\
2 & -0.00040311 & 0 & 0.0020732 & 0 & 0 & 0 & 0 & 0 & 0 & 0 & 0 \\
3 & 0 & -0.0010002 & 0 & 0.0048457 & 0 & 0 & 0 & 0 & 0 & 0 & 0 \\
4 & 0.0061259 & 0 & -0.067032 & 0 & 0.17212 & 0 & 0 & 0 & 0 & 0 & 0 \\
5 & 0 & 0.013607 & 0 & -0.14851 & 0 & 0.38055 & 0 & 0 & 0 & 0 & 0 \\
6 & -0.027586 & 0 & 0.43733 & 0 & -2.2537 & 0 & 3.789 & 0 & 0 & 0 & 0 \\
7 & 0 & -0.069611 & 0 & 1.0873 & 0 & -5.5105 & 0 & 9.1099 & 0 & 0 & 0 \\
8 & 1.4436 & 0 & -31.645 & 0 & 254.75 & 0 & -891.62 & 0 & 1146.2 & 0 & 0 \\
9 & 0 & 3.0921 & 0 & -67.772 & 0 & 545.53 & 0 & -1909.1 & 0 & 2454.1 & 0 \\
10 & -7.2086 & 0 & 191.07 & 0 & -1996.6 & 0 & 10284 & 0 & -26134 & 0 & 26237 \\
\end{tabular}
}
\caption{
CLP $R=8$, Matrix of (rounded) $\alpha_{j,k}$ coefficients, divided by $10^4$
}
\label{tablesix1}
\end{table}

\begin{figure}[htbp!]
\includegraphics[width=\onewidth\textwidth]{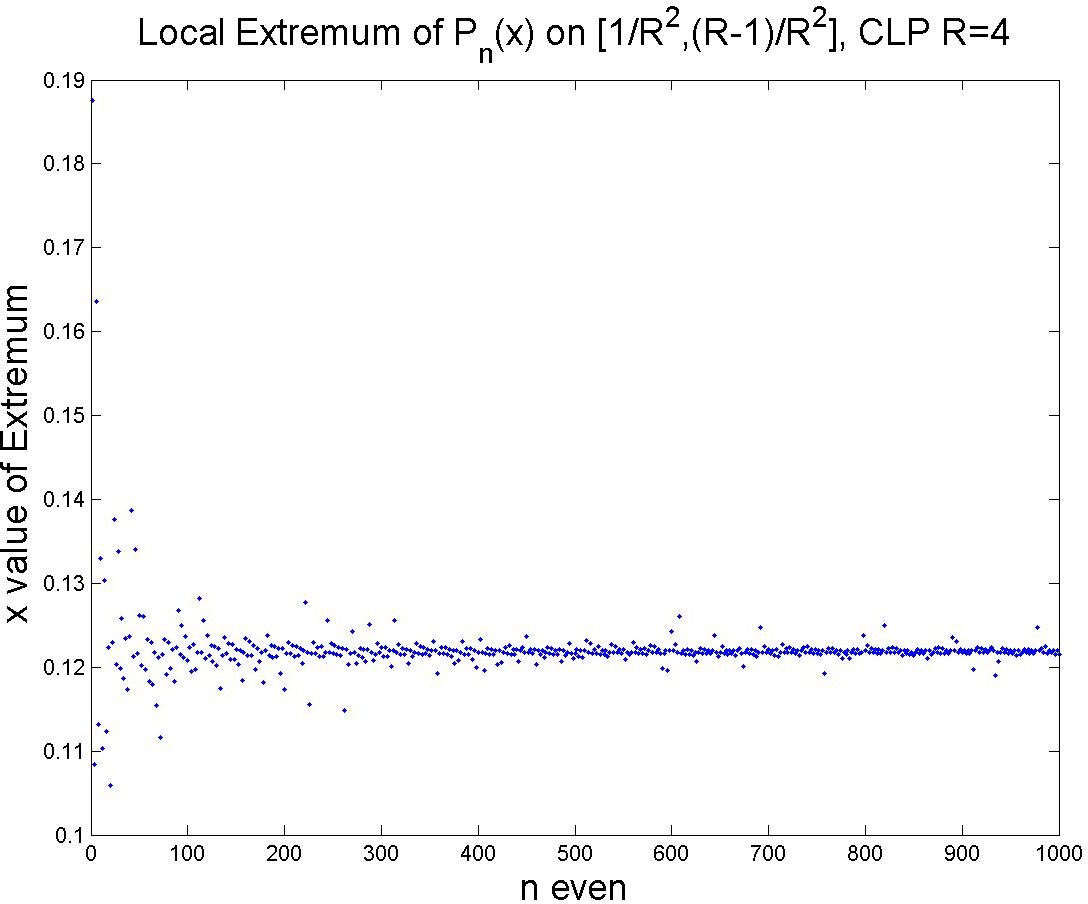}
\caption{$x$-value of local extremum of $P_{n}(x)$ on $[1/R^{2},(R-1)/R^{2}]=[1/16,3/16]$ as a function of $n$, CLP $R=4$.}
\label{figsix2}
\end{figure}

On the other gaps (i.e. the images of $[1/R,(R-1)/R]$ under compositions of the maps $F_{1}$ and $F_{2}$), we observe similar qualitative behavior.  That is, on these smaller gaps, $P_{n}(x)$ resembles a Gaussian for all large and even $n$.  We have seen this self-similar behavior already in Figs. \ref{figthree8} and \ref{figthree9}.  However, the local extremum of $P_{2n}(x)$ on $[1/R^{2},(R-1)/R^{2}]$ is not in general the center of the given gap.  We glimpse this phenomenon in Fig. \ref{figthree10}.  In Fig. \ref{figsix2} we plot the local extremum of $P_{2n}|_{[1/16,3/16]}$ in the CLP $R=4$ case.  The center of this gap occurs at $x=.125$, but the local extremum fluctuates with average around $.122$.

\clpg
\section{WLP as a function of $n$ for fixed $x$}
\label{secwlp}

\begin{figure}[htbp!]
\includegraphics[width=\fourteenwidth\textwidth]{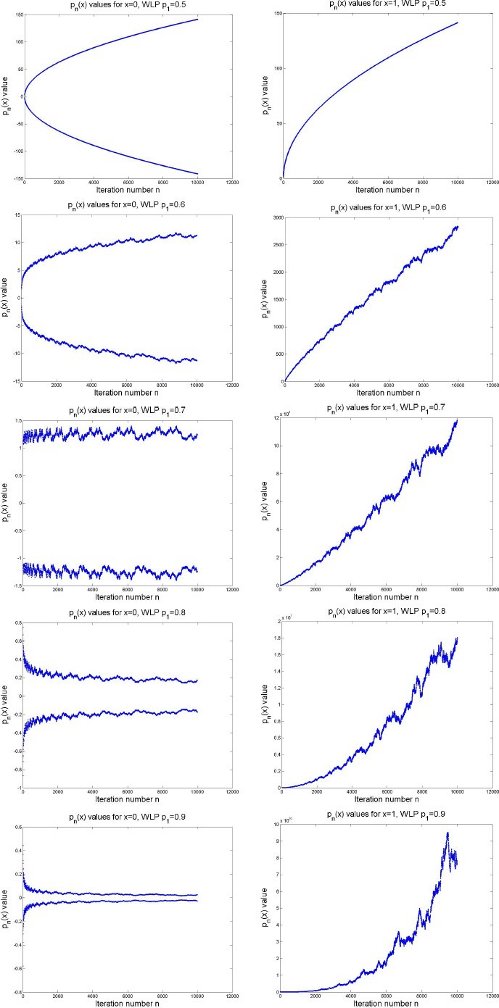}
\caption{$P_{n}(x)$ values vs. $n$ for $x=0,1$ plotted for the WLP $p_{1}=.5,.6,.7,.8,.9$ families}
\label{figseven1}
\end{figure}

\begin{figure}[htbp!]
\includegraphics[width=\onewidth\textwidth]{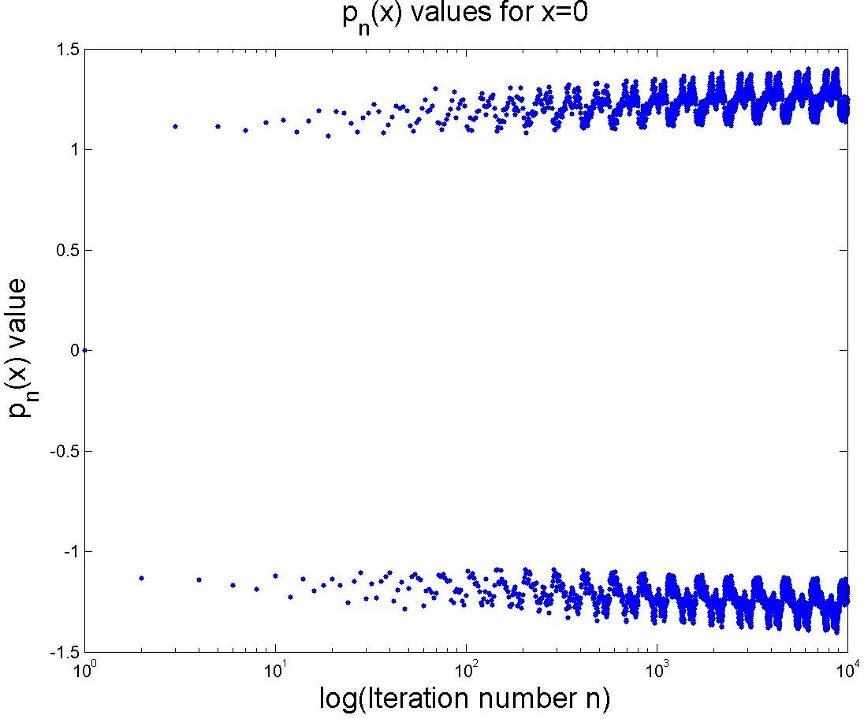}
\caption{$P_{n}(0)$ values vs. $n$ logarithmically scaled in the $x$-axis for WLP $p_{1}=.7$}
\label{figseven1.5}
\end{figure}

For the classical Legendre polynomials on $[0,1]$, we find structure in the values of $P_{n}(x)$ as functions of $n$ for fixed $x$.  As an example, for $x=0$, $P_{n}(x)=(-1)^{n}\sqrt{2n+1}$, and for $x=1$, $P_{n}(1)=\sqrt{2n+1}$.  The computations of $P_{n}(0)$ and of $P_{n}(1)$ are shown at the top of Fig. \ref{figseven1}.  These algebraic relations are a renormalization of the behavior of the classical Legendre polynomials $\widetilde{P}_{n}$ on $[-1,1]$, where we recall that $\widetilde{P}_{n}(1)=1$, $\widetilde{P}_{n}(-1)=(-1)^{n}$ and $\vnormf{\widetilde{P}_{n}}_{L^{2}[-1,1]}^{2}=1/(n+1/2)$ for $n\geq1$ \cite{gautschi04}.  We now increase the WLP weight $p_{1}$ at the points $x=0,1$.  As a result, we see a perturbation of the behavior of the classical Legendre polynomials (on $[0,1]$) in Fig. \ref{figseven1}.  For $p_{1}\approx.7$, a transition seems to occur in the behavior of $\{P_{n}(0)\}_{n=1}^{\infty}$.  At this transition, we appear to have a multiplicative periodic function $P_{n}(0)$, i.e. a function $f\colon\Z^{+}\to\R$ where for some $a>0$ we have $f(ax)=f(x)$ for all $x$.  For $p_{1}=.7$, we can logarithmically scale the $x$-axis as in Fig. \ref{figseven1.5} to see a nearly multiplicative periodic function.  For other $p_{1}$, we evidently have a multiplicative periodic function $P_{n}(0)/n^{\beta}$ for an appropriate choice of $\beta=\beta(p_{1})$.

For generic\footnote{By a generic point we mean a point chosen, with respect to a uniform distribution, among a suitable set of (rational) values in floating point arithmetic.  In all cases, we only treat $x$ representable in double precision floating point arithmetic.} $x\in[0,1]$, if we plot the vectors $(P_{n}(x),P_{n+1}(x))$ in the WLP $p_{1}=.5$ (i.e. classical) case, then these vectors are attracted to an ellipse as $n\to\infty$.  This ellipse is centered at the origin, and its axes and orientation vary with $x$.  Also, the vectors $(P_{n}(x),P_{n+1}(x))$ rotate around the ellipse in a periodic way.  That is, for a certain $k>1$ which depends on $x$ (where $k\in\Z$ is the period), $d((P_{n}(x),P_{n+1}(x)),(P_{n+k}(x),P_{n+k+1}(x)))<\epsilon$ for some small $\epsilon>0$.  Geometrically speaking, $(P_{n}(x),P_{n+1}(x),n)\in\R^{3}$ travels along a helix near the surface of an ellipsoidal cylinder.  Thus, plotting $(P_{n}(x),P_{n+1}(x))$ projects this helix onto the plane.  In Section \ref{secdyn}, we focus on the behavior of $(P_{n}(x),P_{n+1}(x))$ for generic $x\in[0,1]$.  For now, we merely note the contrast between the images of generic and non-generic points.

\begin{figure}[htbp!]
\includegraphics[width=\finonewidth\textwidth]{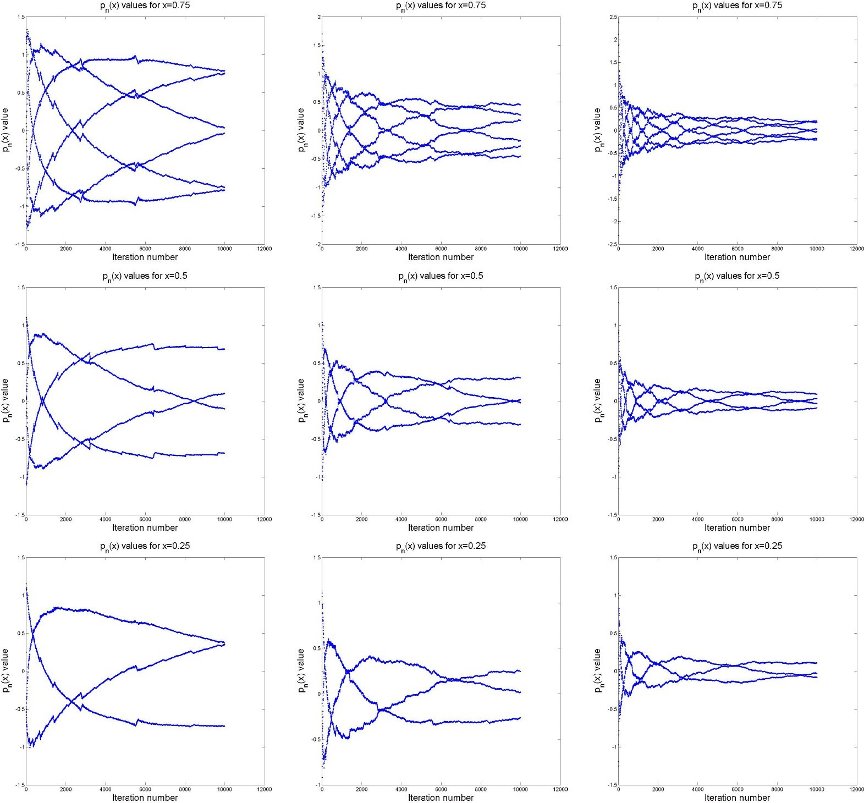}
\caption{$P_{n}(x)$ values vs. $n$ for $x=.75,.5,.25$ plotted for the WLP $p_{1}=.6,.7,.8$ families.  Note that for $p_{1}=.5$, we would get solid horizontal lines.}
\label{figseven2}
\end{figure}

\begin{figure}[htbp!]
\includegraphics[width=\fintwowidth\textwidth]{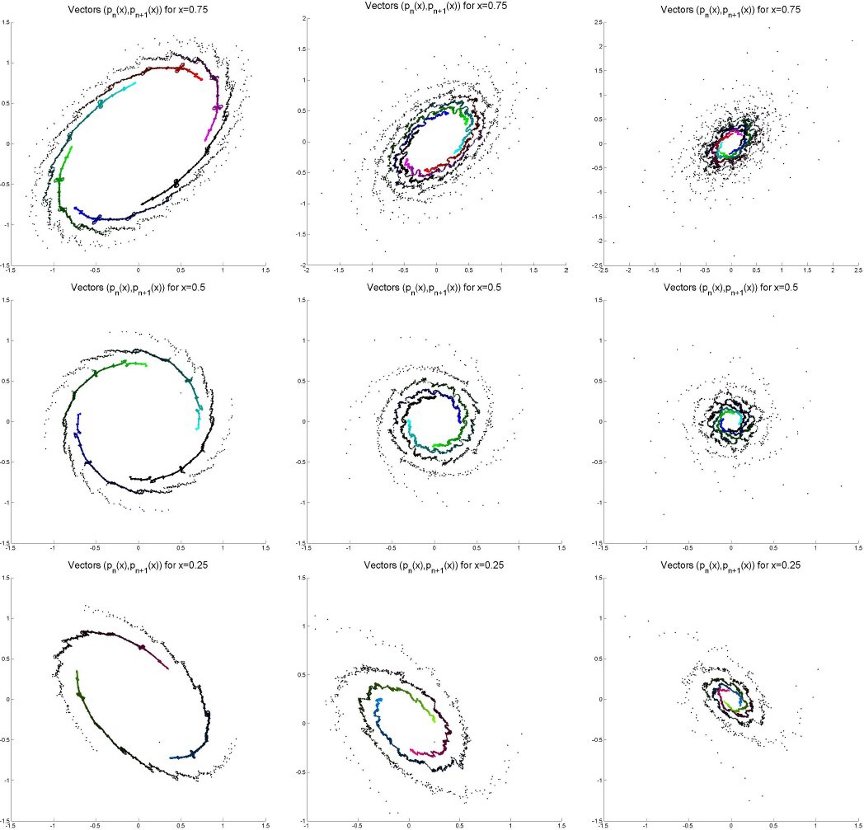}
\caption{Vectors $(P_n(x),P_{n+1}(x))$ for $x=.75,.5,.25$, plotted for the WLP $p_{1}=.6,.7,.8$ families}
\label{figseven3}
\end{figure}

As observed for $x\in\{\{0\},\{1\}\}$ (a set of two non-generic points), increasing $p_{1}$ perturbs the dynamics of the vectors $(P_{n}(x),P_{n+1}(x))$.  Therefore, increasing $p_{1}$ should result in perturbed dynamics for other non-generic $x$.  To this end, consider $x\in\{\{.25\},\{.5\},\{.75\}\}$.  For the classical Legendre polynomials on $[0,1]$, these $x$ values yield finite attractors for $\{(P_{n}(x),P_{n+1}(x))\}_{n=1}^{\infty}$.  These attractors have three, four and six points, respectively.  We call the corresponding number of points $k_{x}$. The vectors $(P_{n}(x),P_{n+1}(x))$ travel in a clockwise fashion about the $k_{x}$ attracting points, as $n$ increases.  When we increase $p_{1}$, the dynamics of $(P_{n}(x),P_{n+1}(x))$ for $x\in\{\{.25\},\{.5\},\{.75\}\}$ change dramatically, as we see in Figs. \ref{figseven2} and \ref{figseven3}.  Instead of cycling ever closer to $k_{x}$ attracting points, the vectors $(P_{n}(x),P_{n+1}(x))$ cycle through $k_{x}$ fractal spiral arms, as we see in Fig. \ref{figseven3}.  However, in the case that $p_{1}>.5$, the attractor of these spiral arms is difficult to determine.  By measuring the distance of $(P_{n}(x),P_{n+1}(x))$ from $(0,0)$, it seems that the spirals are bounded away from the origin as $n\to\infty$.  Therefore, by the symmetry apparent from Figs. \ref{figseven2} and \ref{figseven3}, $(P_{n}(x),P_{n+1}(x))$ should be attracted to a (fractal) ellipse, or a finite set of $k_{x}$ points.  As a further contrast to the $p_{1}=.5$ behavior, we see from Fig. \ref{figseven3} that $(P_{n}(x),P_{n+1}(x),n)\in\R^{3}$ travels along a ``fractal helix'' with radius decreasing in $n$.

\clpg
\section{WLP and CLP Dynamics}
\label{secdyn}

\begin{figure}[htbp!]
\includegraphics[width=\finthreewidth\textwidth]{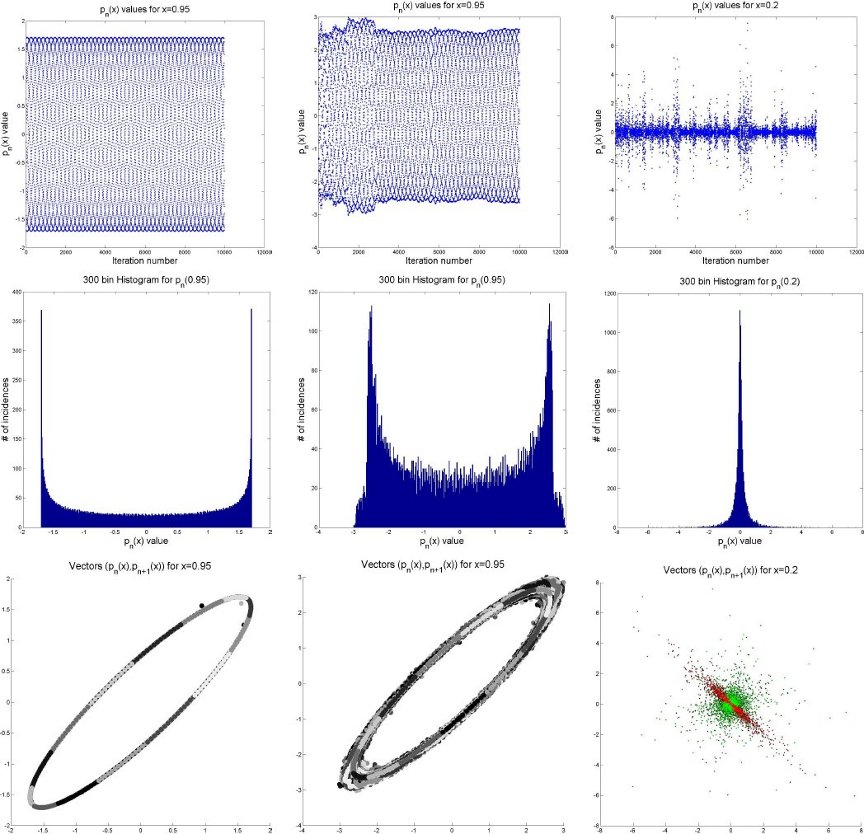}
\caption{Summary of Generic point values: $p_{n}$ values, their histograms, and the vectors $(p_n(x),p_{n+1}(x))$ plotted for the Legendre, WLP $p_{1}=.6$ and CLP, $R=4$  polynomial families}
\label{figeight1}
\end{figure}

As promised, we now examine the dynamics of $P_{n}(x)$ for generic $x$.  Here we examine both WLP and CLP, and we summarize the results in Figure \ref{figeight1}.  For WLP, we recall for $p_{1}=.5$ and fixed, generic $x\in[0,1]$, there exists an integer $k_{x}>1$ so that
\begin{equation}\label{eight1}
d((P_{n}(x),P_{n+1}(x)),(P_{n+k_{x}}(x),P_{n+k_{x}+1}(x)))<\epsilon
\end{equation}
for some small $\epsilon>0$.  Here we mean that the sequence of points\lbreak $\{(P_{n+jk_{x}}(x),P_{n+jk_{x}+1}(x))\}_{j=0}^{\infty}$ is periodic of period $k_{x}$, up to a small error of $\epsilon$ at each step.  We therefore make $k_{x}$ the minimal positive integer satisfying our condition \ref{eight1}. If we plot $P_{n}(x)$ vs. $n$, we find that $P_{n}(x)$ is a superposition of $k_{x}$ sinusoidal functions (see entry $(1,1)$ in Fig. \ref{figeight1}).  These $k_{x}$ (approximately) periodic functions are given by $P_{(k_{x}n+j)}(x)$ for $j\in\{1,2,\ldots,k_{x}\}$.  For example, we can see that $k_{.95}=14$ by counting the number of distinct periodic functions in the plot of $P_{n}(x)$ versus $n$.  Since we essentially have a superposition of phase shifted cosines, it follows that the distribution of values of $P_{n}(x)$ is the function $1/\sqrt{1-y^{2}}$, suitably rescaled (see entry $(2,1)$ in Fig. \ref{figeight1}).  Finally, if we color the iterates $(P_{n}(x),P_{n+1}(x))$ in a $14$-periodic manner, we get entry $(3,1)$ in Fig. \ref{figeight1}.  In this plot, the point $(P_{n+1}(x),P_{n+2}(x))$ overlaps all points of lower index.

So, what happens when we increase $p_{1}$?  As in Section \ref{secwlp}, we observe a perturbation of the $p_{1}=.5$ behavior.  This perturbation becomes more exaggerated the larger $p_{1}$ becomes.  In the middle column of Fig. \ref{figeight1}, we show the same plots described in the previous paragraph, for $p_{1}=.6$.  We again choose $x=.95$ so $k_{.95}=14$, and we observe fourteen vaguely periodic functions.  In other words, for $j\in\{1,2,\ldots,k_{x}\}$, we see that $P_{(k_{x}n+j)}(.95)$ is a perturbed sinusoidal function of $n$.  We can see these perturbations via the distribution function (entry (2,2)) and the plot of $(P_{n}(x),P_{n+1}(x))$ (entry (3,2)) in Fig. \ref{figeight1}.  That is, the distribution of values of $P_{n}(x)$ is a perturbed version of the function $1/\sqrt{1-y^{2}}$, suitably rescaled.  Also, the attractor of the vectors $(P_{n}(x),P_{n+1}(x))$ appears to be a thickened (fractal) ellipse.

The behavior for generic $x$ in the CLP case is displayed in Fig. \ref{figeight1} in the rightmost column.  As an example, we view the iterates of $x=1/5=1/(R+1)\in C_{R}$ for the case $R=4$.  (Recall that $P_{n}(x)$ grows exponentially for $x\in[0,1]\setminus C_{R}$, so we expect more information by viewing $x\in C_{R}$.)  For generic $x$, the distribution of $P_{n}(x)$ values is unclear (Fig. \ref{figeight1}, entry (2,3)).  Nevertheless, the $P_{n}(x)$ values do seem concentrated around the origin.  Also, the vectors $(P_{n}(x),P_{n+1}(x))$ obey a $2$-periodic behavior.  Specifically, the points $\{(P_{2n}(x),P_{2n+1}(x))\}_{n=1}^{\infty}$ cluster along the line $\{(x,y)\in\R^{2}\colon y=-x\}$ (to abuse notation).  Also, the points $\{(P_{2n+1}(x),P_{2n+2}(x))\}_{n=1}^{\infty}$ cluster loosely around the origin.  To observe this behavior, note entry (3,3) in Fig. \ref{figeight1}.

\begin{figure}[htbp!]
\includegraphics[width=\finfourwidth\textwidth]{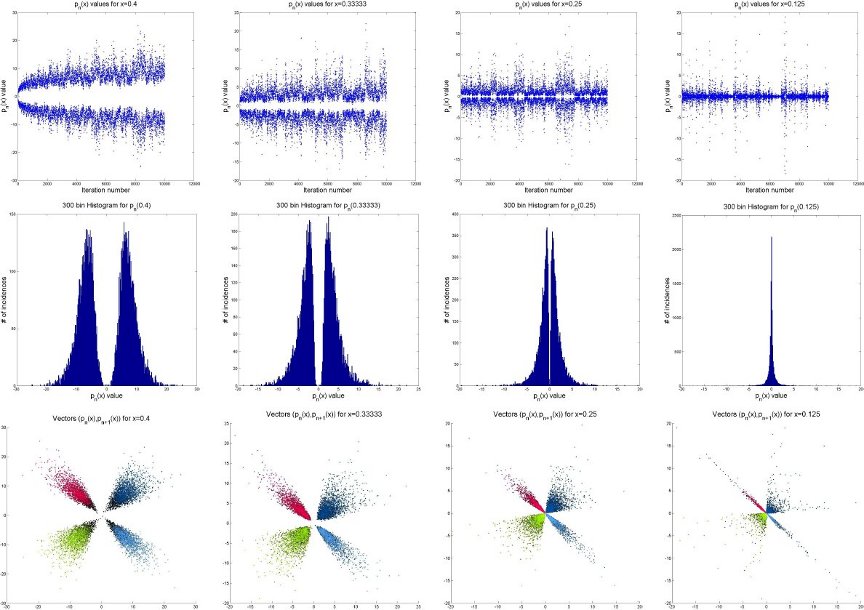}
\caption{$P_{n}$ values, their histograms, and the vectors $(P_n(x),P_{n+1}(x))$ plotted for $x=1/R$, CLP $R=2.5,3,4$ and $8$ polynomial families.}
\label{figeight2}
\end{figure}

In our investigation, we have found only one non-generic point $x\in C_{R}\cap[0,1/2]$ for the CLP case (other than the endpoint $x=0$).  (Note that the symmetry of these $P_{n}(x)$ across $x=1/2$ reduces our investigation to $C_{R}\cap[0,1/2]$.)  Unlike the generic iterates which demonstrate only $2$-periodic behavior, the non-generic point $x=1/R$ exhibits $4$-periodic behavior.  The columns of Fig. \ref{figeight2} correspond to different $p_{1}$ values and the iterates of $x=1/R$.  For $x=1/R$, the vectors $(P_{4n+j}(x),P_{4n+j+1}(x))$ form four disjoint attractors for $j=\{1,2,3,4\}$.  The set $\{(P_{4n+j}(x),P_{4n+j+1}(x))\}_{n=0}^{\infty}$ for $j=1$ resides in the lower left quadrant, for $j=2$ this set lives in the upper left quadrant, for $j=3$ the upper right, and for $j=4$ the lower right.

\clpg
\section{Concluding Discussion}
\label{secconc}

Experimental mathematics allows us to see patterns that otherwise would never be observed.  At times these patterns allow us to formula explicit conjectures, and eventually proofs may be found to turn conjectures into theorems.  But often the story is more complicated.  The observed patterns may turn out, on closer inspection, to be only approximately true.  No clear cut conjectures emerge that are truly supported by the experimental evidence.  In such cases, what should we do?

We could simply discard the experimental results and move on.  But the approximate patterns that we observe may be of great interest.  The experiment might be trying to tell us something that we are not quite able to capture in the conventional format of mathematical statements.  Experimental science is often messy in exactly this way, and we think it would be a shame to limit experimental mathematics to only the clean paradigm.

We see the results reported in this paper in exactly this light.  We do have clean cut conjectures on the $r_{n}$ values in Section \ref{sectwo} and boundedness for WLP in Section \ref{secthree}, but beyond that we only have messy evidence.  We see some CLP Dirichlet kernels that look more like approximate identities than most Dirichlet kernels.  We are able to relate this to the small size of the $r_{n}$ coefficients.  The experimental evidence suggests that we will not find a sequence with $r_{n}$ converging to zero.  So we cannot offer a conjecture analogous to the results of \cite{strichartz06} about uniform convergence of certain special partial sums of the orthogonal polynomial expansion of an arbitrary continuous function.  However, the results of \cite{strichartz06} suggest that perhaps such a statement might hold only for special values of $R$, and we have only examined a few choices of $R$.

The approximate equalities for CLP discussed in Section \ref{secapprox} appear very striking when one looks at the graphs of the polynomials.  It is on closer inspection that one observes the deviation from equality.  A clean conjecture might say that the deviation goes to zero in some limit, but the evidence does not really support such a conjecture.  Nevertheless, we find the approximate equalities intrinsically interesting.

Similarly, a glance at the graphs of CLP on the gaps suggests they are close to being Gaussian.  Of course one comes to expect Gaussian limits in many mathematical situations.  But here in Section \ref{secclp}, on closer inspection, we see a decided deviation from the Gaussian model.  In this case we believe there is a true limit law; we just have not been able to find it.

Perhaps the most interesting contribution of this paper is the dynamical perspective discussed in Sections \ref{secwlp} and \ref{secdyn}.  Here we may be excused from offering explicit conjectures on the grounds that the pictures offer a view of significant structures that were previously unrecognized.  We hope that others will be inspired to investigate the dynamics perspective in a systematic way.

\clpg

\end{document}